\documentclass[12pt,oneside,reqno]{amsart}
\usepackage{graphicx}
\usepackage{mathrsfs}
\usepackage{stmaryrd}
\usepackage{amsfonts}
\usepackage{cite}
\usepackage{enumerate,amsmath,amssymb,amsthm}
\usepackage{booktabs} 
\usepackage{diagbox} 
\usepackage{xcolor}

\newcommand{\dif}{\mathrm{d}}

\newcommand{\be}{\begin{eqnarray}}
	\newcommand{\ee}{\end{eqnarray}}
\newcommand{\ce}{\begin{eqnarray*}}
	\newcommand{\de}{\end{eqnarray*}}
\newtheorem{theorem}{Theorem}[section]
\newtheorem{lemma}[theorem]{Lemma}
\newtheorem{remark}[theorem]{Remark}
\newtheorem{definition}[theorem]{Definition}
\newtheorem{proposition}[theorem]{Proposition}
\newtheorem{Examples}[theorem]{Examples}
\newtheorem{corollary}[theorem]{Corollary}
\newtheorem{condition}[theorem]{Condition}

\def\t{\theta}

\def\d{\delta}

\def\g{\gamma}
\def\s{\sigma}
\def\eps{\varepsilon}

\def\[{{\Big[}}
\def\]{{\Big]}}
\def\<{{\langle}}
\def\>{{\rangle}}
\def\({{\Big(}}
\def\){{\Big)}}

\def\W{{\mathcal W}}

\def\no{\nonumber}
\def\bt{\begin{theorem}}
	\def\et{\end{theorem}}
\def\bl{\begin{lemma}}
	\def\el{\end{lemma}}
\def\br{\begin{remark}}
	\def\er{\end{remark}}
\def\bx{\begin{Examples}}
	\def\ex{\end{Examples}}
\def\bd{\begin{definition}}
	\def\ed{\end{definition}}
\def\bp{\begin{proposition}}
	\def\ep{\end{proposition}}
\def\bc{\begin{corollary}}
	\def\ec{\end{corollary}}
\def\bco{\begin{condition}}
	\def\eco{\end{condition}}

\def\cC{{\mathcal C}}
\def\cD{{\mathcal D}}

\def\cG{{\mathcal G}}

\def\cN{{\mathcal N}}

\def\mE{{\mathbb E}}

\def\mK{{\mathbb K}}

\def\mN{{\mathbb N}}

\def\mP{{\mathbb P}}

\def\mR{{\mathbb R}}

\def\mW{{\mathbb W}}

\def\sA{{\mathscr A}}
\def\sB{{\mathscr B}}
\def\sC{{\mathscr C}}

\def\sF{{\mathscr F}}

\def\sL{{\mathscr L}}

\def\sP{{\mathscr P}}

\def\sV{{\mathscr V}}

\def\geq{\geqslant}
\def\leq{\leqslant}

\begin{document}

\allowdisplaybreaks
\title{Well-posedness for path-dependent multivalued McKean-Vlasov stochastic differential equations}

\author{Ying Ma and Huijie Qiao}

\thanks{{\it AMS Subject Classification(2020):} 60H10}

\thanks{{\it Keywords:} Path-dependent multivalued McKean-Vlasov stochastic differential equations, maximal monotone operators, non-Lipschitz conditions.}

\thanks{This work was supported by NSF of China (No.12071071) and the Jiangsu Provincial Scientific Research Center of Applied Mathematics (No. BK20233002).}

\thanks{Corresponding author: Huijie Qiao, hjqiaogean@seu.edu.cn}

\subjclass{}

\date{}

\dedicatory{Department of Mathematics,
Southeast University,\\
Nanjing, Jiangsu 211189, P.R.China}

\begin{abstract}
This work concerns a type of path-dependent multivalued McKean-Vlasov stochastic differential equations. First of all, we prove the well-posedness for path-dependent multivalued stochastic differential equations under the Lipschitz conditions. Then by constructing Lipschitz approximation sequences, we generalize the result to the case of the non-Lipschitz conditions. Finally, based on the obtained results, the well-posedness for path-dependent multivalued McKean-Vlasov stochastic differential equations under the non-Lipschitz conditions is established by iterating in distributions.
\end{abstract}

\maketitle \rm

\section{Introduction}
McKean-Vlasov stochastic differential equations (SDEs for short), also known as \\distribution-dependent SDEs or mean-field SDEs, were initiated by McKean \cite{McKean} who was inspired by Kac's program in Kinetic Theory. These equations are closely associated with nonlinear Fokker-Planck equations (cf. \cite{hrockw}). And there exist extensive results regarding McKean-Vlasov SDEs, for instance \cite{dq1,dq2,hrenw, lq, wangf1} and the references therein. McKean-Vlasov SDEs with maximal monotone operators are called multivalued McKean-Vlasov SDEs. At present, there have been some results about multivalued McKean-Vlasov SDEs, such as the well-posedness (\cite{chi, qg}), ergodicity (\cite{q}) and large deviation principles (\cite{flqz}).

Meanwhile, there is a type of SDEs whose coefficients depend on the past path. These SDEs are called path-dependent SDEs, alternatively referred to as functional SDEs or SDEs with memory. As for path-dependent SDEs, the well-posedness have been established under global Lipschitz conditions (cf. \cite{Moha}) or under local Lipschitz conditions and linear growth conditions (cf. \cite{Mao}). Later, von Renesse and Scheutzow \cite{Mvon} weakened the Lipschitz conditions and showed that only weak one-sided Lipschitz conditions are sufficient for the well-posedness. Wang \cite{wangf2} also obtained the well-posedness under the log-Lipschitz conditions. Recently, for special path-dependent SDEs, Huang \cite{huangx2} proved the well-posedness under the singular conditions.

Next, let $(\Omega,\sF,\{\sF_t\}_{t\geq 0},\mP)$ be a complete filtered probability space and $\{W(t),t\geq0\}$ be a $m$-dimensional standard Brownian motion on it. Fix $T>0$ and consider the following multivalued SDE:
\be\left\{\begin{array}{ll}
	X(t)\in -A(X(t))\dif t+b(t,X_t,\sL_{X_t})\dif t+\s(t,X_t,\sL_{X_t})\dif W(t), 0\leq t\leq T,\\
	X_0=\xi \in \sC,
\end{array}
\right.
\label{1eq}
\ee
where $A:\mathbb{R}^d\to 2^{\mathbb{R}^d}$ is a maximal monotone operator (See Subsection \ref{maxi}), $b:[0,T]\times\sC\times\sP_{2}^{\sC}\to \mR^d$, $ \sigma:[0,T]\times\sC\times\sP_{2}^{\sC}\to \mR^{d\times m}$ are Borel measurable and $\xi$ is a $\sF_0$-measurable $\sC$-valued random variable. Here $X_{\cdot}$ is the segment functional of $X(\cdot)$ defined as $X_t(\theta)=X(t+\theta),\theta\in[-r_0,0]$ and $\sL_{X_t}$ is the distribution of $X_t$. Since the coefficients of Eq.(\ref{1eq}) depend on not only paths but also distributions, we call Eq.(\ref{1eq}) a path-dependent multivalued McKean-Vlasov SDE. Moreover, these equations like Eq.(\ref{1eq}) have appeared in biology, mechanics, epidemiology, and finance which are inevitably memory-dependent (cf. \cite{wxz}).

If $A=0$, Eq.(\ref{1eq}) becomes into the following path-dependent McKean-Vlasov SDE
\be
\dif X(t)=b(t,X_t,\sL_{X_t})\dif t+\sigma(t,X_t,\sL_{X_t})\dif W(t), \quad t\geq 0.
\label{2eq}
\ee
For Eq.(\ref{2eq}), Huang, R\"ockner and Wang \cite{hrockw} firstly obtained the well-posedness under one-sided Lipschitz conditions by applying a distribution iteration technique. Later, if the diffusion coefficient $\sigma(t,X_t,\sL_{X_t})$ is replaced by $\sigma(t,X(t),\sL_{X_t})$, Huang \cite{huangx1} and Zhao \cite{zhaox} established the well-posedness under the singular conditions.

If $A\neq0$, there are few results about Eq.(\ref{1eq}). Therefore, the aim of this paper is to establish the well-posedness for Eq.(\ref{1eq}) under non-Lipschitz conditions. To do so, we first obtain the well-posedness for the path-dependent multivalued SDEs under the Lipschitz conditions. Then by constructing Lipschitz approximation sequences, we generalize the result to the case of the non-Lipschitz conditions. Our result can cover \cite[Theorem 4.1.1]{wangf2} and \cite[Theorem 2.3]{Mvon} partly. Finally, based on the obtained results, the well-posedness for the path-dependent multivalued McKean-Vlasov SDEs is proved under the non-Lipschitz conditions by iterating in distribution. We mention that even although $A=0$, our result can {\it not} be covered by \cite{huangx1} and \cite{zhaox}. Besides, we emphasis that the appearance of the maximal monotone operator $A$ renders the derivation intricate.

The rest of the paper is organized as follows. In Section \ref{preliminaries}, we recall some notations, concepts and some useful lemmas. After that, we prove the well-posedness for path-dependent multivalued SDEs under Lipschitz conditions and non-Lipschitz conditions respectively in Section \ref{pdexiuni}. Finally in Section \ref{pddexiuni}, the well-posedness for path-dependent multivalued McKean-Vlasov SDEs under the non-Lipschitz conditions is established.

The following convention will be used throughout the paper: $C$ with or without indices will denote different positive constants whose values may change from one place to another.

\section{Preliminaries}\label{preliminaries}

In this section, we recall some notations, concepts and some useful lemmas.

\subsection{Notations}
We shall use $|\cdot|$ and $\|\cdot\|$ for norms of vectors and matrices, respectively. Let $\<\cdot,\cdot\>$ denote scalar product in $\mathbb{R}^d$. 

Fix $r_0>0$. For any closed subset ${\rm D}\subset\mR^d$, let $\sC:=C([-r_0,0]; {\rm D})$ be the collection of all the continuous functions from $[-r_0,0]$ to ${\rm D}$ with the uniform norm $\|\xi\|_\infty:=\sup\limits_{\theta\in[-r_0,0]}|\xi(\theta)|$. For $\gamma(\cdot)\in C([-r_0,\infty);\rm D)$, the segment functional $\gamma_\cdot\in C(\mathbb{R}_{+};\sC)$ is defined as $\gamma_t(\theta):=\gamma(t+\theta),\ \theta\in[-r_0,0],t\geq 0.$

Let $\sP_{2}^{\sC}$ be the class of probability measures on $\sC$ which have finite second-order moments, that is 
$$
\mu(\|\cdot\|_{\infty}^{2}):=\int_{\sC}\|\xi\|_{\infty}^2\mu(\dif \xi)<\infty, \quad \mbox{for}~ \mu\in\sP_{2}^{\sC}.
$$
 Then $\sP_{2}^{\sC}$ is a polish space under the Wasserstein distance
$$
\mW_2\left(\mu,\nu\right):=\inf_{\pi\in\cC(\mu,\nu)}\left(\int_{\sC\times\sC}\|\xi-\eta\|_\infty^2\pi(\dif\xi,\dif\eta)\right)^{\frac{1}{2}},
$$
where $\cC(\mu,\nu)$ is the class of couplings for $\mu$ and $\nu$. Moreover, if $\mu$ and $\nu$ are distributions of $X,Y\in L^2(\Omega, \sF, \mP; \sC)$ respectively, it holds that 
\ce
\mathbb{W}_2(\mu,\nu)^2\leq \mE\|X-Y\|^2_{\infty},
\de
where $\mE$ denotes the expectation with respect to $\mP$.

\subsection{Maximal monotone operators}\label{maxi}
For a multivalued operator $A:\mathbb{R}^d\to 2^{\mathbb{R}^d}$, let
\ce
\mathcal{D}(A)&:=&\{x\in\mathbb{R}^d:A(x)\neq\emptyset\},\\
\mathrm{Gr}(A)&:=&\{(x,y)\in\mathbb{R}^{2d}:x\in\mathcal{D}(A),y\in A(x)\}.
\de
We say that $A$ is monotone if for any $(x_{1},y_{1}),(x_{2},y_{2})\in \mathrm{Gr}(A)$, 
$$
\< x_{1}-x_{2},y_{1}-y_{2}\>\geq 0;
$$
$A$ is maximal monotone if  
$$
(x_1,y_1)\in \mathrm{Gr}(A) \Longleftrightarrow\<x_1-x_2,y_1-y_2\>\geq 0,\forall(x_2,y_2)\in \mathrm{Gr}(A).
$$

For any $T>0$, let $\sV_0$ be the set of all continuous functions $K:[0,T]\to \mathbb{R}^d$ with finite variations and $K(0)=0$. And we denote the variation of $K$ on $[0,s]$ for some $s\in[0,T]$ as $|K|_0^{s}$. Set
\ce	
&&\sA:=\{\left(X,K\right):X\in C([0,T];\overline{\cD(A)}),K\in\sV_{0},\\ 
&&\qquad\quad \text{and}~\< X(t)-x,\dif K(t)-y\dif t\>\geq 0,\forall(x,y)\in\mathrm{Gr}(A)\}.
\de

As for $\sA$, the following properties hold (cf. \cite{Cepa98,zhangx2}):
\bl\label{L1}
For $X\in C([0,T];\overline{\cD(A)})$ and $K\in\sV_0$, the following statements are equivalent:
\begin{enumerate}[(i)]
\item $(X,K)\in\sA$;
\item For any $(x,y)\in C([0,T];\mathbb{R}^d)$ with $(x(t),y(t))\in\mathrm{Gr}(A)$, it holds that 
$$
\< X(t)-x(t),\dif K(t)-y(t)\dif t\>\geq 0;
$$
\item For any $(X',K')\in\sA$, it holds that 
$$
\< X(t)-X'(t),\dif K(t)-\dif K'(t)\>\geq 0.
$$
\end{enumerate}
\el

\bl\label{L2}
Assume that $\mathrm{Int}(\cD(A))\not=\emptyset$, then for any $a\in\mathrm{Int}(\cD(A))$, there exists constants $\gamma_1>0$ and $\gamma_2,\gamma_3\geq0$ such that for any $(X,K)\in\sA$ and $0\leq s<t\leq T$,
$$
\int_s^t\< X(r)-a,\dif K(r)\>\geq\gamma_1|K|_s^t-\gamma_2\int_s^t\left|X(r)-a\right|\dif r-\gamma_3(t-s).
$$	
\el

\bl\label{L3}
Assume that $\{K^n,n\in\mathbb{N}\}\subset \sV_0$ converge to some $K\in C([0,T];\mathbb{R}^d)$ and $\sup\limits_{n\in\mathbb{N}}|K^{n}|_{0}^{T}<\infty$. Then $K\in\sV_0$ and 
$$
\lim_{n\to\infty}\int_0^T\<X^n(s),\mathrm{d}K^n(s)\>=\int_0^T\<X(s),\mathrm{d}K(s)\>,
$$
where $\{X^n\}\subset C([0,T];\mathbb{R}^d)$ converges to some $X\in C([0,T];\mathbb{R}^d)$.
\el

\section{The well-posedness of path-dependent multivalued SDEs}\label{pdexiuni}

In this section, we take ${\rm D}=\overline{\cD(A)}$, i.e. $\sC=C([-r_0,0], \overline{\cD(A)})$ and study the well-posedness of path-dependent multivalued SDEs.

For any $T>0$, let $(\Omega,\sF,\{\sF_t\}_{t\in[0,T]},\mP)$ be a complete filtered probability space and $\{W(t),t\in[0,T]\}$ be a $m$-dimensional standard Brownian motion on it. Consider the following path-dependent multivalued SDE:
\be\left\{\begin{array}{ll}
\dif X(t)\in -A(X(t))\dif t+f(t,X_t)\dif t+g(t,X_t)\dif W(t),\quad t\in[0,T],\\
X_0=\xi \in \sC,
\end{array}
\right.
\label{eq1}
\ee
where $A:\mathbb{R}^d\to 2^{\mathbb{R}^d}$ is a maximal monotone operator with $\mathrm{Int}(\cD(A))\not=\emptyset$, $f: [0,T]\times\sC\to \mR^d, g: [0,T]\times\sC\to \mR^{d\times m}$ are Borel measurable, $X_{\cdot}$ is the segment functional of $X(\cdot)$ and $\xi$ is a $\sF_0$-measurable $\sC$-valued random variable.

We first give the definitions of strong solutions, weak solutions and the pathwise uniqueness for Eq.(\ref{eq1}).

\bd(Strong solutions)\label{strosolu}
We say that Eq.(\ref{eq1}) admits a strong solution with the initial value $\xi$ if there exists a pair of adapted processes $(X_\cdot,K(\cdot))$ on the filtered probability space  $(\Omega,\sF,\{\sF_t\}_{t\in[0,T]},\mP)$ such that
\begin{enumerate}[(i)]
\item $\mP(X_0=\xi)=1$,
\item $X(t):=X_t(0)\in\sF_t^{W}$, where $\sF_t^{W}$ is the $\sigma$-field filtration generated by $W$, 
\item $(X_\cdot(0),K(\cdot))\in\sA\ a.s.\ \mP$,
\item it holds that
$$
\mP\left\{\int_0^T(|f(s,X_s)|+\|g(s,X_s)\|^2)\dif s<\infty\right\}=1
$$
and
$$
X_t(0)=\xi(0)-K(t)+\int_0^t f(s,X_s)\dif s+\int_0^tg(s,X_s)\dif W(s), \quad t\in[0,T].
$$
\end{enumerate}
\ed

Next, let $\sL_{\xi}$ be the distribution of $\xi$.

\bd(Weak solutions)\label{weaksolu}
We say that  Eq.(\ref{eq1}) admits a weak solution with the initial law $\sL_{\xi}\in \sP(\mR^{d})$, if there exists a filtered probability space $(\hat{\Omega}, \hat{\sF}, \{\hat{\sF}_t\}_{t\in[0,T]}, \hat{\mP})$, a $m$-dimensional standard $(\hat{\sF}_{t})$-Brownian motion $\hat{W}$ as well as a pair of $(\hat{\sF}_{t})$-adapted process $(\hat{X}_\cdot,\hat{K}(\cdot))$ defined on $(\hat{\Omega}, \hat{\sF}, \{\hat{\sF}_t\}_{t\in[0,T]}, \hat{\mP})$ such that
\begin{enumerate}[(i)]
	\item $\hat{\mP}\circ\hat{X}^{-1}_0=\sL_{\xi}$,
	\item $(\hat{X}_\cdot(0),\hat{K}(\cdot))\in \sA$ a.s. $\hat{\mP}$,
	\item it holds that  
	$$
\mP\left\{\int_0^T(\mid{f(s,\hat{X}_{s})}\mid+\parallel{g(s,\hat{X}_{s})}\parallel^2)\dif s<+\infty\right\}=1,
	$$
	and
	$$
	\hat{X}_t(0)=\hat{X}_0(0)-\hat{K}(t)+\int_{0}^{t}f(s,\hat{X}_{s})\dif s +\int_{0}^{t}g(s,\hat{X}_{s})\dif \hat{W}_{s},\quad 0\leq{t}\leq{T}.
	$$
\end{enumerate}
\ed

Such a solution will be denoted  by $(\hat{\Omega}, \hat{\sF}, \{\hat{\sF}_t\}_{t\in[0,T]}, \hat{\mP}; \hat{W}, (\hat{X}_\cdot,\hat{K}(\cdot)))$.

\bd(Pathwise Uniqueness) \label{pathuniq}
Suppose $(\hat{\Omega}, \hat{\sF}, \{\hat{\sF}_t\}_{t\in[0,T]}, \hat{\mP}; \hat{W}, (\hat{X}^1_\cdot,\hat{K}^1(\cdot)))$ and $(\hat{\Omega}, \hat{\sF}, \{\hat{\sF}_t\}_{t\in[0,T]}, \hat{\mP}; \hat{W}, (\hat{X}^2_\cdot,\hat{K}^2(\cdot)))$ are two weak solutions for Eq.(\ref{eq1}) with $\hat{X}^1_{0}=\hat{X}^2_0$. If $\hat{\mP}\{(\hat{X}^1_t,\hat{K}^1(t))=(\hat{X}^2_t,\hat{K}^2(t)), t\in[0,T]\}=1$, we say that the pathwise uniqueness holds for Eq.(\ref{eq1}).
\ed

\subsection{Well-posedness of path-dependent multivalued SDEs under the Lipschitz conditions}
In this subsection, we prove the well-posedness of path-dependent multivalued SDEs under the Lipschitz conditions.

Assume:
\begin{enumerate}[$({\bf H}_1)$] 
\item There exists a constant $L_1>0$ such that for any $t\in[0,T]$ and any $\zeta\in\sC$
$$
|f(t,\zeta)|^2+\|g(t,\zeta)\|^2\leq L_1(1+\|\zeta\|^2_{\infty}).
$$
\item There exists a constant $L_2>0$ such that for any $t\in[0,T]$ and any $\zeta,\eta\in\sC$
$$
|f(t,\zeta)-f(t,\eta)|^2+\|g(t,\zeta)-g(t,\eta)\|^2\leq L_2\|\zeta-\eta\|^2_\infty.
$$
\end{enumerate}

Now, we state the main result in this subsection.

\bt\label{wellpose1}
Suppose that $({\bf H}_1)$ and $({\bf H}_2)$ hold. Then for any $\sF_0$-measurable initial value $X_0=\xi$ satisfying $\mE\|\xi\|_\infty^2<\infty$, Eq.(\ref{eq1}) has a unique strong solution $(X_\cdot,K(\cdot))$ and $\mE\(\sup\limits_{-r_0 \leq t \leq T}|X(t)|^2\)<\infty.$
\et

In order to prove the above theorem, we construct an iteration sequence. Set $X^{(0)}_t(\theta)=\xi(0\wedge(t+\theta))$ for $  \theta\in[-r_0,0], t\in[0,T]$. For $n\in{\mN^+}$, consider the following multivalued SDE:
\be\left\{\begin{array}{ll}
\dif X^{(n)}(t)\in -A(X^{(n)}(t))\dif t + f(t,X_{t}^{(n-1)})\dif t + g(t,X_{t}^{(n-1)})\dif W(t),\quad t\in[0,T],\\
X_{0}^{(n)}= \xi\in\sC.
\end{array}
\right.
\label{aux}
\ee
Note that Eq.(\ref{aux}) is a mutivalued SDE whose coefficients do not depend on its path. Thus, by \cite{Cepa95}, Eq.(\ref{aux}) has a unique strong solution denoted as $(X^{(n)}_\cdot,K^{(n)}(\cdot))$. We intend to prove that $(X^{(n)}_\cdot, K^{(n)}(\cdot))$ converges as $n\rightarrow\infty$ and its limit is a weak solution of Eq.(\ref{eq1}). As for $(X^{(n)}_\cdot, K^{(n)}(\cdot))$, we have the following lemma.

\bl\label{property-0}
Assume $({\bf H}_1)$ and $({\bf H}_2)$, it holds that

$(i)$
\be 
\sup\limits_{n\in{\mN^+}}\mE\(\sup\limits_{s\in[-r_0,t]}|X^{(n)}(s)|^2\)+\sup\limits_{n\in{\mN^+}}\mE|K^{(n)}|_0^t<\infty, \quad t\in[0,T];          
\label{bounded}
\ee

$(ii)$ there exists a $t_0>0$ such that for any $n\in{\mN^+}$
\be
\mE\sup_{t\in[-r_0,t_0]}|X^{(n+1)}(t)-X^{(n)}(t)|^2\leq2(\frac{1}{2})^{n}\left(\mE\sup_{t\in[-r_0,t_0]}|X^{(1)}(t)|^2+\mE\|\xi\|^2_\infty\right);
\label{iteration}
\ee
 
$(iii)$ $\{X^{(n)}_{\cdot}\}$ is a Cauchy sequence in $L^2(\Omega, \sF, \mP; C([0,t_0],\sC))$.
\el
\begin{proof}
We first prove $(i)$. For $n=1$ and any $a\in \mathrm{Int}(\cD(A))$, by the It\^o formula, Lemma \ref{L2} and $({\bf H}_1)$, it holds that
\ce
&&\left|X^{(1)}(t)-a\right|^2\\
&=&\left|\xi(0)-a\right|^{2}-2\int_{0}^{t}\< X^{(1)}(s)-a,\dif K^{(1)}(s)\> +2\int_{0}^{t}\< X^{(1)}(s)-a,f(s,X_{s}^{(0)})\>\dif s\\
&&+2\int_{0}^{t}\< X^{(1)}(s)-a,g(s,X_{s}^{(0)})\dif W(s)\>+\int_{0}^{t}\left\|g(s,X_{s}^{(0)})\right\|^{2}\dif s\\
&\leq&|\xi(0)-a|^{2}-2\gamma_{1}|K^{(1)}|_{0}^{t}+2\gamma_{2}\int_{0}^{t}|X^{(1)}(s)-a|\dif s+2\gamma_{3}t+\int_0^t|X^{(1)}(s)-a|^2\dif s\\
&&+\int_0^t|f(s,X_{s}^{(0)})| ^2\dif s+\int_0^t\left\|g(s,X_{s}^{(0)})\right\|^2\dif s+2\int_{0}^{t}\< X^{(1)}(s)-a,g(s,X_{s}^{(0)})\dif W(s)\>\\
&\leq&|\xi(0)-a|^2-2\gamma_1|K^{(1)}|_0^t+2\gamma_2\int_0^t\left(1+|X^{(1)}(s)-a|^2\right)\dif s+2\gamma_3t +\int_{0}^{t}|X^{(1)}(s)-a|^{2}\dif s\\
&&+L_{1}\int_{0}^{t}\left(1+\|X^{(0)}_{s}\|_{\infty}^{2}\right)\dif s+2\int_0^t\< X^{(1)}(s)-a,g(s,X_{s}^{(0)})\dif W(s)\> \\
&=&|\xi(0)-a|^{2}-2\gamma_{1}|K^{(1)}|_{0}^{t}+2\gamma_{2}t+2\gamma_{3}t+\left(2\gamma_{2}+1\right)\int_{0}^{t}|X^{(1)}(s)-a|^{2}\dif s \\
&&+L_{1}\int_{0}^{t}\left(1+\|X^{(0)}_{s}\|_{\infty}^{2}\right)\dif s+2\int_{0}^{t}\< X^{(1)}(s)-a,g(s,X_{s}^{(0)})\dif W(s)\>.
\de	
The Burkholder-Davis-Gundy inequality and the Young inequality imply that
\ce
&& \mE\sup_{s\in[0,t]}\left|X^{(1)}(s)-a\right|^2+2\gamma_1\mE\left|K^{(1)}\right|_0^t \\
&\leq&\mE\left|\xi(0)-a\right|^2+2\gamma_{2}t+2\gamma_{3}t+\left(2\gamma_{2}+1\right)\int_{0}^{t}\mE|X^{(1)}(s)-a|^{2}\dif s \\
&&+L_{1}\int_{0}^{t}\mE\left(1+\|X^{(0)}_{s}\|_{\infty}^{2}\right)\dif s +12\mE\(\int_0^t\left|X^{(1)}(s)-a\right|^2\left\|g(s,X_s^{(0)})\right\|^2\dif s\)^{1/2} \\
&\leq&\mE\left|\xi(0)-a\right|^2+2T(\gamma_{2}+\gamma_{3})+\left(2\gamma_{2}+1\right)\int_{0}^{t}\mE|X^{(1)}(s)-a|^{2}\dif s\\
&&+L_{1}\int_{0}^{t}\mE\left(1+\|X^{(0)}_{s}\|_{\infty}^{2}\right)\dif s+\frac{1}{2}\mE\(\sup_{s\in[0,t]}\left|X^{(1)}(s)-a\right|^2\) \\
&&+C\mE\int_0^t\left\|g(s,X_s^{(0)})\right\|^2\dif s.
\de
Then by $({\bf H}_1)$ we get that
\be
&& \mE\sup_{s\in[0,t]}\left|X^{(1)}(s)-a\right|^2+4\gamma_1\mE\left|K^{(1)}\right|_0^t \no\\
&\leq&2\mE\left|\xi(0)-a\right|^2+4T(\gamma_{2}+\gamma_{3})+2\left(2\gamma_{2}+1\right)\int_{0}^{t}\mE\sup_{u\in[-r_0,s]}\left|X^{(1)}(u)-a\right|^2\dif s \no\\
&&+2L_1(C+1)\mE\int_0^t\left(1+\left\|X_s^{(0)}\right\|_\infty^2\right)\dif s\no\\
&\leq&2\mE\left|\xi(0)-a\right|^2+4T(\gamma_{2}+\gamma_{3})+2L_1(C+1)T(1+2|a|^2) \no\\
&&+2\left(2\gamma_{2}+1\right)\int_{0}^{t}\mE\sup_{u\in[-r_0,s]}\left|X^{(1)}(u)-a\right|^2\dif s+4L_1(C+1)\int_{0}^{t}\mE\sup_{u\in[-r_0,s]}\left|X^{(0)}(u)-a\right|^2\dif s\no\\
&\leq&2\mE\left|\xi(0)-a\right|^2+2T\left[2(\gamma_{2}+\gamma_{3})+L_1(C+1)(1+2|a|^2)\right] \no\\
&&+2\left(2\gamma_{2}+1\right)\int_{0}^{t}\mE\sup_{u\in[-r_0,s]}\left|X^{(1)}(u)-a\right|^2\dif s+4L_1(C+1)t\mE\sup_{u\in[-r_0,t]}\left|X^{(0)}(u)-a\right|^2.
\label{x1k1x0}
\ee

Besides, it holds that
\ce
\sup_{s\in[-r_0,t]}\left|X^{(1)}(s)-a\right|^2\leq \sup_{s\in[-r_0,0]}\left|X^{(1)}(s)-a\right|^2+\sup_{s\in[0,t]}\left|X^{(1)}(s)-a\right|^2.
\de
Since $X^{(1)}(s)=X^{(1)}_0(s)=\xi(s)$ for $s\in[-r_0,0]$, we infer that
\ce
\sup_{s\in[-r_0,t]}\left|X^{(1)}(s)-a\right|^2\leq \sup_{s\in[-r_0,0]}\left|\xi(s)-a\right|^2+\sup_{s\in[0,t]}\left|X^{(1)}(s)-a\right|^2,
\de
which together with (\ref{x1k1x0}) implies that
\be
&&\mE\sup_{s\in[-r_0,t]}\left|X^{(1)}(s)-a\right|^2+4\gamma_1\mE\left|K^{(1)}\right|_0^t\no\\
 &\leq&\mE\sup_{s\in[-r_0,0]}\left|\xi(s)-a\right|^2+2\mE\left|\xi(0)-a\right|^2+2T\left[2(\gamma_{2}+\gamma_{3})+L_1(C+1)(1+2|a|^2)\right] \no\\
&&+2\left(2\gamma_{2}+1\right)\int_{0}^{t}\mE\sup_{u\in[-r_0,s]}\left|X^{(1)}(u)-a\right|^2\dif s\no\\
&&+4L_1(C+1)t\mE\sup_{u\in[-r_0,t]}\left|X^{(0)}(u)-a\right|^2.
\label{n=1}
\ee
Note that $X^{(0)}_t(\theta)=\xi(0\wedge(t+\theta))$. Thus, it holds that
\ce
\mE\sup_{u\in[-r_0,t]}\left|X^{(0)}(u)-a\right|^2=\mE\sup_{u\in[-r_0,0]}\left|\xi(u)-a\right|^2.
\de
The Gronwall inequality yields that
\ce
&&\mE\sup\limits_{u\in[-r_0,t]}|X^{(1)}(u)-a|^2\\
&\leq& \left[\(3+4L_1(C+1)T\)\mE\sup_{u\in[-r_0,0]}\left|\xi(u)-a\right|^2+C_1\right]e^{C_2t}\\
&\leq&  \left[2\(3+4L_1(C+1)T\)\mE\|\xi\|^2_\infty+2\(3+4L_1(C+1)T\)|a|^2+C_1\right]e^{C_2t},
\de
where 
\ce
C_{1}&=&2T\left[2(\gamma_{2}+\gamma_{3})+L_1(C+1)(1+2|a|^2)\right],\\
C_{2}&=&2\left(2\gamma_{2}+1\right).
\de
By this and (\ref{n=1}) we conclude that
\ce
\mE\sup_{s\in[-r_0,t]}|X^{(1)}(s)|^{2}+\mE\left|K^{(1)}\right|_0^t<\infty.
\de

Now assume that (\ref{bounded}) holds for $n=k$ for some $k\geq1$ and we prove that it holds for $n=k+1$. By the similar derivation to that for (\ref{n=1}), for any $a\in \mathrm{Int}(\cD(A))$, we have that
\ce
&&\mE\sup_{s\in[-r_0,t]}\left|X^{(k+1)}(s)-a\right|^2+4\gamma_1\mE\left|K^{(k+1)}\right|_0^t \\
&\leq&\mE\sup_{s\in[-r_0,0]}\left|\xi(s)-a\right|^2+2\mE\left|\xi(0)-a\right|^2+2T\left[2(\gamma_{2}+\gamma_{3})+L_1(C+1)(1+2|a|^2)\right]\no\\
&&+2\left(2\gamma_{2}+1\right)\int_{0}^{t}\mE\sup_{u\in[-r_0,s]}\left|X^{(k+1)}(u)-a\right|^2\dif s +4L_1(C+1)t\mE\sup_{u\in[-r_0,t]}\left|X^{(k)}(u)-a\right|^2,
\de
which together with the Gronwall inequality leads to
$$
\mE\sup_{s\in[-r_0,t]}|X^{(k+1)}(s)|^{2}+\mE\left|K^{(k+1)}\right|_0^t <\infty.
$$
Thus, by the mathematical induction method we conclude that (\ref{bounded}) holds for all $n\in\mN^+$.

Next, we prove $(ii)$. The It\^o formula, Lemma \ref{L1} and $({\bf H}_2)$ imply that
\ce
&&\left|X^{(n+1)}(t)-X^{(n)}(t)\right|^2\\ 
&=&-2\int_0^t\< X^{(n+1)}(s)-X^{(n)}(s),\dif K^{(n+1)}(s)-\dif K^{(n)}(s)\> \\
&&+2\int_0^t\< X^{(n+1)}(s)-X^{(n)}(s),f(s,X_s^{(n)})-f(s,X_s^{(n-1)})\>\dif s \\
&&+2\int_0^t\< X^{(n+1)}(s)-X^{(n)}(s),(g(s,X_s^{(n)})-g(s,X_s^{(n-1)}))\dif W(s)\> \\
&&+\int_0^t\|g(s,X_s^{(n)})-g(s,X_s^{(n-1)})\|^2\dif s \\
&\leq&\int_{0}^{t}|X^{(n+1)}(s)-X^{(n)}(s)|^{2}\dif s+\int_{0}^{t}|f(s,X_{s}^{(n)})-f(s,X_s^{(n-1)})|^{2}\dif s
\\
&&+\int_0^t\left\|g(s,X_s^{(n)})-g(s,X_s^{(n-1)})\right\|^2\dif s\\
&&+2\int_0^t\< X^{(n+1)}(s)-X^{(n)}(s),(g(s,X_s^{(n)})-g(s,X_s^{(n-1)}))\dif W(s)\> \\
&\leq&\int_{0}^{t}|X^{(n+1)}(s)-X^{(n)}(s)|^{2}\dif s+L_2\int_{0}^{t}\|X^{(n)}_s-X^{(n-1)}_s\|^2_{\infty}\dif s\\
&&+2\int_0^t\< X^{(n+1)}(s)-X^{(n)}(s),(g(s,X_s^{(n)})-g(s,X_s^{(n-1)}))\dif W(s)\>.
\de
It follows from the Burkholder-Davis-Gundy inequality, the Young inequality and $({\bf H}_2)$ that
\be
&&\mE\sup_{s\in[0,t]}|X^{(n+1)}(s)-X^{(n)}(s)|^{2} \no\\ 
&\leq&\int_0^t\mE|X^{(n+1)}(s)-X^{(n)}(s)|^2\dif s+L_2\int_{0}^{t}\mE\sup_{u\in[-r_0,s]}|X^{(n)}(u)-X^{(n-1)}(u)|^2\dif s\no\\
&&+12\mE\left(\int_{0}^{t}|X^{(n+1)}(s)-X^{(n)}(s)|^2\|g(s,X^{(n)}_s)-g(s,X^{(n-1)}_s)\|^2\dif s\right)^{1/2}\no\\
&\leq&\int_0^t\mE\sup_{u\in[-r_0,s]}|X^{(n+1)}(u)-X^{(n)}(u)|^2\dif s+L_2\int_{0}^{t}\mE\sup_{u\in[-r_0,s]}|X^{(n)}(u)-X^{(n-1)}(u)|^2\dif s\no\\
&&+\frac{1}{2}\mE\sup_{s\in[0,t]}|X^{(n+1)}(s)-X^{(n)}(s)|^2+C\mE\int_0^t\|g(s,X_s^{(n)})-g(s,X_s^{(n-1)})\|^2\dif s\no\\
&\leq&\int_0^t\mE\sup_{u\in[-r_0,s]}|X^{(n+1)}(u)-X^{(n)}(u)|^2\dif s+(L_2+CL_2)\int_{0}^{t}\mE\sup_{u\in[-r_0,s]}|X^{(n)}(u)-X^{(n-1)}(u)|^2\dif s\no\\
&&+\frac{1}{2}\mE\sup_{s\in[0,t]}|X^{(n+1)}(s)-X^{(n)}(s)|^2.
\label{n1n}
\ee
Then by some calculation we know that
\be
&&\mE\sup_{s\in[0,t]}|X^{(n+1)}(s)-X^{(n)}(s)|^{2}\no\\
&\leq&2\int_0^t\mE\sup_{u\in[-r_0,s]}|X^{(n+1)}(u)-X^{(n)}(u)|^2\dif s+2(L_2+CL_2)\int_{0}^{t}\mE\sup_{u\in[-r_0,s]}|X^{(n)}(u)-X^{(n-1)}(u)|^2\dif s\no\\
&\leq&2\int_0^t\mE\sup_{u\in[-r_0,s]}|X^{(n+1)}(u)-X^{(n)}(u)|^2\dif s+2(L_2+CL_2)t\mE\sup_{s\in[-r_0,t]}|X^{(n)}(s)-X^{(n-1)}(s)|^2.
\label{n+1nn-1}
\ee

Besides, it holds that
\ce
\sup_{s\in[-r_0,t]}|X^{(n+1)}(s)-X^{(n)}(s)|^{2}\leq \sup_{s\in[-r_0,0]}|X^{(n+1)}(s)-X^{(n)}(s)|^{2}+\sup_{s\in[0,t]}|X^{(n+1)}(s)-X^{(n)}(s)|^{2}.
\de
Note that 
\ce
X^{(n+1)}(s)=X^{(n)}(s)=\xi(s), \quad s\in[-r_0,0].
\de
Thus, we infer that
\ce
\sup_{s\in[-r_0,t]}|X^{(n+1)}(s)-X^{(n)}(s)|^{2}\leq\sup_{s\in[0,t]}|X^{(n+1)}(s)-X^{(n)}(s)|^{2},
\de
which together with (\ref{n+1nn-1}) and the Gronwall inequality yields that
\ce
\mE\sup_{s\in[-r_0,t]}|X^{(n+1)}(s)-X^{(n)}(s)|^{2}\leq2(L_2+CL_2)te^{2t}\mE\sup_{s\in[-r_0,t]}|X^{(n)}(s)-X^{(n-1)}(s)|^{2}.
\de
Taking $t_0>0$ such that $2(L_2+CL_2)t_0e^{2t_0}\leq\frac{1}{2}$, we conclude that
\ce
\mE\sup_{s\in[-r_0,t_0]}|X^{(n+1)}(s)-X^{(n)}(s)|^{2}&\leq&\frac{1}{2}\mE\sup_{s\in[-r_0,t_0]}|X^{(n)}(s)-X^{(n-1)}(s)|^{2}\\
&\leq&(\frac{1}{2})^n\mE\sup_{s\in[-r_0,t_0]}|X^{(1)}(s)-X^{(0)}(s)|^{2}\\
&\leq&2(\frac{1}{2})^n\left(\mE\sup_{s\in[-r_0,t_0]}|X^{(1)}(s)|^2+\mE\|\xi\|^{2}_{\infty}\right).
\de
So, (\ref{iteration}) is proved.

Finally, we prove $(iii)$. By (\ref{iteration}), for any $n,m\in\mN^+$ and $m>n$, it holds that
\ce
\left(\mE\sup_{s\in[-r_0,t_{0}]}|X^{(n)}(s)-X^{(m)}(s)|^{2}\right)^{1/2}
&\leq&\sum_{k=n}^{m-1}\left(\mE\sup_{s\in[-r_0,t_{0}]}\left|X^{k+1}(s)-X^{k}(s)\right|^{2}\right)^{1/2} \\
&\leq&\sum_{k=n}^{m-1}2^{1/2}2^{-k/2}\left(\mE\sup_{s\in[-r_0,t_{0}]}|X^{(1)}(s)|^{2}+\mE\|\xi\|^2_\infty\right)^{1/2}.
\de
Thus, we obtain that
\ce
\lim\limits_{n,m\rightarrow\infty}\mE\sup_{s\in[-r_0,t_{0}]}|X^{(n)}(s)-X^{(m)}(s)|^{2}=0,
\de
which implies that $\{X^{(n)}_{\cdot}\}$ is a Cauchy sequence in $L^2(\Omega, \sF, \mP; C([0,t_0],\sC))$. The proof is complete.

Now, it is the position to prove Theorem \ref{wellpose1}.

{\bf Proof of Theorem \ref{wellpose1}.} We first prove the existence. According to $(iii)$ in Lemma \ref{property-0}, there exists a process $(X_{t})_{t\in[0,t_0]}\in L^2(\Omega, \sF, \mP; C([0,t_0],\sC))$ with $X_0=\xi$ satisfying 
$$
\lim_{n\to\infty}\mE\sup\limits_{t\in[0,t_0]}\left\|X^{(n)}_t-X_t\right\|^2_{\infty}=0.
$$
Moreover, by $({\bf H}_1)$,$({\bf H}_2)$ and the dominated convergence theorem, we can get that
\ce
\int_{0}^{t}f(s,X_s^{(n)})\dif s &\to&\int_{0}^{t}f(s,X_{s})\dif s, \\
\int_{0}^{t}g(s,X_s^{(n)})\dif W(s)&\to&\int_{0}^{t}g(s,X_{s})\dif W(s).
\de
Let 
\ce
K(t):=\xi(0)+\int_0^t f(s,X_s)\dif s+\int_0^tg(s,X_s)\dif W(s)-X(t).
\de
Note that 
\ce
K^{(n)}(t)=\xi(0)+\int_0^t f(s,X^{(n-1)}_s)\dif s+\int_0^tg(s,X^{(n-1)}_s)\dif W(s)-X^{(n)}(t).
\de 
Thus, it follows from (\ref{bounded}) and Lemma \ref{L3} that $K\in \sV_0$ and $(X(\cdot), K(\cdot))\in\sA$. Thus $(X_\cdot,K(\cdot))$ is a weak solution of Eq.(\ref{eq1}) up to time $t_0$. By solving the equation piecewise in time, we conclude that Eq.(\ref{eq1}) has a weak solution on $[0,T]$ with
\ce
\mE\sup_{s\in[-r_0,T]}|X(s)|^2<\infty.
\de  

Finally, we prove the pathwise uniqueness for Eq.(\ref{eq1}). Assume that $(\Omega, \sF, \{\sF_t\}_{t\in[0,T]}, \mP; \\W, (Y_\cdot,\tilde{K}(\cdot)))$ is the other weak solution for Eq.(\ref{eq1}) with the same initial value $Y_0=\xi$, that is,
\ce
Y(t)=\xi(0)-\tilde{K}(t)+\int_0^t f(s,Y_s)\dif s+\int_0^tg(s,Y_s)\dif W(s).
\de
By the similar deduction to that for (\ref{n1n}), we have that
\ce
&&\mE\(\sup_{s\in[0,t]}|X(s)-Y(s)|^2\) \\
&\leq&(1+L_2+CL_2)\int_0^t\mE\left(\sup_{u\in[-r_0,s]}|X(u)-Y(u)|^2\right)\dif s+\frac{1}{2}\mE\left(\sup_{s\in[0,t]}|X(s)-Y(s)|^2\right),
\de 
and
\ce
\mE\(\sup_{s\in[0,t]}|X(s)-Y(s)|^{2}\)\leq2(1+L_2+CL_2)\int_0^t\mE\(\sup_{u\in[-r_0,s]}|X(u)-Y(u)|^2\)\dif s.
\de
Besides, it holds that
\ce
X(s)=Y(s)=\xi, \quad s\in[-r_0,0],
\de
which yields that
\ce
\sup_{s\in[-r_0,t]}|X(s)-Y(s)|^{2}=\sup_{s\in[0,t]}|X(s)-Y(s)|^{2}.
\de
Therefore, we infer that
\ce
\mE\(\sup_{s\in[-r_0,t]}|X(s)-Y(s)|^{2}\)\leq2(1+L_2+CL_2)\int_0^t\mE\(\sup_{u\in[-r_0,s]}|X(u)-Y(u)|^2\)\dif s.
\de
Then the Gronwall inequality implies that $\mE\(\sup\limits_{s\in[-r_0,t]}\left|X(s)-Y(s)\right|^2\)=0, \forall t\in[0,T]$. That is, $X(t)=Y(t),\forall t\in[-r_0,T]\ a.s.\ \mathbb{P}$.

For $\forall t\in[0,T]$,
\ce
K(t) &=&\xi(0)-X(t)+\int_{0}^{t}f(s,X_{s})\dif s+\int_{0}^{t}g(s,X_{s})\dif W(s)\no \\
&=&\xi(0)-Y(t)+\int_{0}^{t}f(s,Y_{s})\dif s+\int_{0}^{t}g(s,Y_{s})\dif W(s)\no \\
&=&\tilde{K}(t).
\de
Moreover, $K(t), \tilde{K}(t)$ are continuous in $t$. Thus, $K(t)=\tilde{K}(t),\forall t\in[0,T]\  a.s.\mathbb{P}$. The proof is complete.
\end{proof}

\subsection{Well-posedness of path-dependent multivalued SDEs under the non-Lipschitz conditions}

In this subsection, we study the well-posedness of path-dependent multivalued SDEs under the non-Lipschitz conditions.

Assume:
\begin{enumerate}[$({\bf H}'_1)$] 
\item There exists a constant $L'_1>0$ such that for any $t\in[0,T]$ and any $\zeta\in\sC$
$$
|f(t,\zeta)|^2+\|g(t,\zeta)\|^2\leq L'_1(1+\|\zeta\|^2_{\infty}).
$$ 
\item For $\forall t\in[0,T]$, $f(t,\cdot)$ is continuous on $\sC$ and there exists a constant $L'_2>0$ such that for any $t\in[0,T]$ and any $\zeta,\eta\in\sC$
\ce
&&2\< f(t,\zeta)-f(t,\eta),\zeta(0)-\eta(0)\>\leq L'_2\kappa(\|\zeta-\eta\|^2_\infty),\\
&&\|g(t,\zeta)-g(t,\eta)\|^2\leq L'_2\kappa(\|\zeta-\eta\|^2_\infty),
\de
where $\kappa(x)$ is a positive, strictly increasing, continuous concave function and satisfies
\ce
\kappa(0)=0, \quad \int_{0+}\frac{1}{\kappa(x)}\dif x=\infty.
\de
\end{enumerate}

At present, we formulate the main result in this subsection.

\bt\label{wellpose2}
Suppose that $({\bf H}'_1)$ and $({\bf H}'_2)$ hold. Then for any $\sF_0$ measurable initial value $X_0=\xi$ satisfying $\mE\|\xi\|_\infty^2<\infty$, Eq.(\ref{eq1}) has a unique strong solution $(X_\cdot,K(\cdot))$ with $\mE\(\sup\limits_{-r_0 \leq t \leq T}|X(t)|^2\)<\infty$.
\et

In order to prove the above theorem, we prepare a key proposition.

\bp\label{bsboun}
Assume that $f, g$ are bounded and satisfy $({\bf H}'_2)$. Then for any $\sF_0$ measurable initial value $X_0=\xi$ satisfying $\mE\|\xi\|_\infty^2<\infty$, Eq.(\ref{eq1}) has a unique strong solution $(X_\cdot,K(\cdot))$ with $\mE\(\sup\limits_{-r_0 \leq t \leq T}|X(t)|^2\)<\infty$.
\ep
\begin{proof}
The proof is divided into four steps. In the first step we construct an approximation sequence of Eq.(\ref{eq1}) and prove the tightness of their solutions in the second step. The existence and pathwise uniqueness of weak solutions for Eq.(\ref{eq1}) are established in third and fourth steps, respectively.

{\bf Step 1.} We construct an approximation sequence of Eq.(\ref{eq1}).

For any $n\geq 1$, set
\ce
\varphi_n(\zeta)(s):=n\int_s^{s+\frac{1}{n}}\frac{I_{[-r_0,1]}(r)(\|\zeta\|_\infty\wedge n)}{\|\zeta\|_\infty}\zeta(r\wedge 0)\dif r, \quad \zeta\in \sC, \quad s\in[-r_0,0].
\de
Define 
\ce
&&f_n(t, \zeta):=\mE f\left(t,\varphi_n(\zeta)+\frac{1}{n}\tilde{W}(r_0+\cdot)\right),\\
&&g_n(t, \zeta):=\mE g\left(t,\varphi_n(\zeta)+\frac{1}{n}\tilde{W}(r_0+\cdot)\right).
\de
where $\tilde{W}$ is a $d$-dimensional Brownian motion with $\tilde{W}(0)=0$. So, for every $n\geq 1$, $f_n, g_n$ satisfy
\ce
&&|f_n(t, \zeta)|^2+\|g_n(t, \zeta)\|^2\leq C, \quad t\in[0,T], \zeta\in \sC\\
&&|f_n(t, \zeta)-f_n(t, \eta)|^2+\|g_n(t, \zeta)-g_n(t, \eta)\|^2\leq C_n\|\zeta-\eta\|^2_\infty, \quad t\in[0,T], \zeta, \eta\in \sC,
\de
where $C$ is independent of $n$ and $C_n$ depends on $n$. By Theorem \ref{wellpose1}, the following path-dependent multivalued SDE
\ce\left\{\begin{array}{ll}
\dif X^{(n)}(t)\in -A(X^{(n)}(t))\dif t+f_n(t,X^{(n)}_t)\dif t+g_n(t,X^{(n)}_t)\dif W(t),\quad t\in[0,T],\\
X^{(n)}_0=\xi \in \sC,
\end{array}
\right.
\de
has a unique strong solution $(X^{(n)}_{\cdot}, K^{(n)}(\cdot))$. 

{\bf Step 2.} We prove that $\{(X^{(n)}_\cdot, K^{(n)}(\cdot)), n\in\mN\}$ is tight in $C([0,T],\sC)\times C([0,T],\mR^d)$.

In order to obtain that $\{(X^{(n)}_\cdot, K^{(n)}(\cdot)), n\in\mN\}$ is tight in $C([0,T],\sC)\times C([0,T],\mR^d)$, we only need to deduce that

$(i)$
 \ce
\sup\limits_{n\in{\mN^+}}\mE\(\sup\limits_{s\in[-r_0,T]}|X^{(n)}(s)|^2\)+\sup\limits_{n\in{\mN^+}}\mE|K^{(n)}|_0^T\leq C;
\de

$(ii)$ for any $\t>0$
\be
&&\lim\limits_{\d\downarrow 0}\sup\limits_{n}\mP\left\{\sup\limits_{-r_0\leq s\leq t\leq s+\d\leq T}|X^{(n)}(t)-X^{(n)}(s)|\geq\t\right\}=0,\label{pxeuxes}\\
&&\lim\limits_{\d\downarrow 0}\sup\limits_{n}\mP\left\{\sup\limits_{0\leq s\leq t\leq s+\d\leq T}|K^{(n)}(t)-K^{(n)}(s)|\geq\t\right\}=0. \label{pkeukes}
\ee

First of all, by the same deduction to that for (\ref{bounded}), we conclude $(i)$. Then, we prove $(ii)$. By the It\^o formula, we have that for $0\leq s\leq t\leq s+\d\leq T$
\be
&&\left|X^{(n)}(t)-X^{(n)}(s)\right|^2 \no\\
&=&-2\int_{s}^{t}\< X^{(n)}(u)-X^{(n)}(s),\dif K^{(n)}(u)\>+ 2\int_{s}^{t}\< X^{(n)}(u)-X^{(n)}(s),f_n(u,X^{(n)}_u)\>\dif u\no\\
&& +2\int_{s}^{t}\< X^{(n)}(u)-X^{(n)}(s),g_n(u,X^{(n)}_u)\dif W(u)\> +\int_{s}^{t}\left\|g_n(u,X^{(n)}_u)\right\|^{2}\dif u\no\\
&\leq&-2\int_{s}^{t}\< X^{(n)}(u)-X^{(n)}(s),\dif K^{(n)}(u)\>+\int_{s}^{t}\left|X^{(n)}(u)-X^{(n)}(s)\right|^{2}\dif s+2C\d \no\\
&&+2\int_{s}^{t}\< X^{(n)}(u)-X^{(n)}(s),g_n(u,X^{(n)}_u)\dif W(u)\>.
\label{xntskes}
\ee

In the following, we estimate $-2\int_{s}^{t}\< X^{(n)}(u)-X^{(n)}(s),\dif K^{(n)}(u)\>$. For $\eps>0$ and $R>0$, let
$$
K_{\eps,R}:=\{x\in\mR^d:\forall y\notin\overline{\cD(A)},|x-y|\geq\eps\ \text{and}\ |x-a_0|\leq R\},
$$
where $a_0\in\mathrm{Int}(\cD(A))$ such that $K_{\eps,R}\not=\emptyset$ for all $R>0$ and all $\eps\in(0,\eps_0)$, for a certain $\eps_0>0$ that does not depend on $R$. Set 
$$
f_R(\eps):=\sup\{|x^{*}|:x^{*}\in A(x), x \in K_{\eps,R}\},
$$
and $|f_R(\eps)|<\infty$ since $A$ is locally bounded on $\mathrm{Int}(\cD(A))$. Again set for any $\delta>0$, 
$$
g_R(\delta):=\inf\{\eps\in(0,\eps_0):f_R(\eps)\leq\delta^{-1/2}\},
$$
and $f_R(\delta+g_R(\delta))\leq \delta^{-1/2}$ and $\lim\limits_{\delta\downarrow0}g_R(\delta)=0$.

Let $\delta_R>0$ such that $\delta_R+g_R(\delta_R)<\eps_0$. For fixed $R>0$ and $\delta\in(0,\delta_R\wedge 1]$, $K_{\delta+g_R(\delta),R}\not=\emptyset$. Let $X^{n,\delta,R}(s)$ denote the projection of $X^{(n)}(s)$ on $K_{\delta+g_R(\delta),R}$. Therefore, on the set $\{\sup\limits_{t\in[0,T]}|X^{n}(t)-a_0|\leq R\}$, it holds that
$$
|X^{n,\delta,R}(s)-X^{(n)}(s)|\leq \delta+g_R(\delta).
$$
Hence, on the set $\{\sup\limits_{t\in[0,T]}|X^{n}(t)-a_0|\leq R\}$, if $0\leq t-s\leq \delta$, we have that
\be
&&-2\int_{s}^{t}\<X^{(n)}(u)-X^{(n)}(s),\dif K^{(n)}(u)\>\no\\
&=&-2\int_{s}^{t}\<X^{(n)}(u)-X^{n,\delta,R}(s),\dif K^{(n)}(u)\>-2\int_{s}^{t}\<X^{n,\delta,R}(s)-X^{(n)}(s),\dif K^{(n)}(u)\>\no\\
&\leq&-2\int_{s}^{t}\<X^{(n)}(u)-X^{n,\delta,R}(s),z^{*}\>\dif u+2(\delta+g_R(\delta))|K^{(n)}|_0^T\no\\
&\leq&4(R+|a_0|)\delta^{1/2}+2(\delta+g_R(\delta))|K^{(n)}|_0^T,
\label{dK}
\ee 
where $z^{*}\in A(X^{n,\delta,R}(s))$. Moreover, by (\ref{xntskes}) and (\ref{dK}), we infer that 
\be
&&\sup\limits_{s\leq t\leq s+\d}\left|X^{(n)}(t)-X^{(n)}(s)\right|^2I_{\{\sup\limits_{t\in[0,T]}|X^{n}(t)-a_0|\leq R\}} \no\\
&\leq&4(R+|a_0|)\delta^{1/2}+2(\delta+g_R(\delta))|K^{(n)}|_0^T+4R^2\d+2C\d \no\\
&&+2(|a_0|+R)\sup\limits_{s\leq t\leq s+\d}\left|\int_{s}^{t}g_n(u,X^{(n)}_u)\dif W(u)\right|\no\\
&&+2\sup\limits_{s\leq t\leq s+\d}\left|\int_{s}^{t}\< X^{(n)}(u),g_n(u,X^{(n)}_u)\dif W(u)\>\right|.
\label{stsdxnr}
\ee

In the following, by the Chebyshev inequality it holds that
\ce
&&\mP\left\{\sup\limits_{0\leq s\leq t\leq s+\d\leq T}|X^{(n)}(t)-X^{(n)}(s)|\geq\t\right\}\\
&=&\mP\left\{\sup\limits_{0\leq s\leq t\leq s+\d\leq T}|X^{(n)}(t)-X^{(n)}(s)|\geq\t,\sup\limits_{t\in[0,T]}|X^{n}(t)-a_0|\leq R\right\}\\
&&+\mP\left\{\sup\limits_{0\leq s\leq t\leq s+\d\leq T}|X^{(n)}(t)-X^{(n)}(s)|\geq\t,\sup\limits_{t\in[0,T]}|X^{n}(t)-a_0|>R\right\}\\
&\leq&\mP\left\{\sup\limits_{0\leq s\leq t\leq s+\d\leq T}\left|X^{(n)}(t)-X^{(n)}(s)\right|^2I_{\{\sup\limits_{t\in[0,T]}|X^{n}(t)-a_0|\leq R\}} \geq \t^2\right\}\\
&&+\frac{2(C+|a_0|^{2})}{R^2}.
\de
Now, we estimate the first term in the right hand of the above inequality. (\ref{stsdxnr}) and the Chebyshev inequality imply that 
\ce
&&\mP\left\{\sup\limits_{0\leq s\leq t\leq s+\d\leq T}\left|X^{(n)}(t)-X^{(n)}(s)\right|^2I_{\{\sup\limits_{t\in[0,T]}|X^{n}(t)-a_0|\leq R\}} \geq \t^2\right\}\\
&\leq&\mP\left\{4(R+|a_0|)\delta^{1/2}+2(\delta+g_R(\delta))|K^{(n)}|_0^T+4R^2\d+2C\d\geq \t^2/3\right\}\\
&&+\mP\left\{2(|a_0|+R)\sup\limits_{s\leq t\leq s+\d}\left|\int_{s}^{t}g_n(u,X^{(n)}_u)\dif W(u)\right|\geq \t^2/3\right\}\\
&&+\mP\left\{2\sup\limits_{s\leq t\leq s+\d}\left|\int_{s}^{t}\< X^{(n)}(u),g_n(u,X^{(n)}_u)\dif W(u)\>\right|\geq \t^2/3\right\}\\
&\leq&\frac{3}{\t^2}\left(4(R+|a_0|)\delta^{1/2}+2(\delta+g_R(\delta))C+4R^2\d+2C\d\right)\\
&&+\mP\left\{2(|a_0|+R)\sup\limits_{s\leq t\leq s+\d}\left|\int_{s}^{t}g_n(u,X^{(n)}_u)\dif W(u)\right|\geq \t^2/3\right\}\\
&&+\mP\left\{2\sup\limits_{s\leq t\leq s+\d}\left|\int_{s}^{t}\< X^{(n)}(u),g_n(u,X^{(n)}_u)\dif W(u)\>\right|\geq \t^2/3\right\}.
\de
Besides, note that for any $n\geq 1$, $\(\int_{s}^{t}g_n(u,X^{(n)}_u)\dif W(u)\)_{s\leq t\leq s+\d}$ is a local martingale, and
$$
\int_{s}^{s}g_n(u,X^{(n)}_u)\dif W(u)=0.
$$
Moreover, by the boundedness of $g$, we know that $\left\{\(\int_s^\cdot\|g_n(u,X^{(n)}_u)\|^2\dif u\), n\in\mN\right\}$ is tight. So, \cite[Lemma 1]{zhengw} guarantees that $\left\{\(\int_s^\cdot g_n(u,X^{(n)}_u)\dif W(u)\), n\in\mN\right\}$ is tight and 
$$
\lim\limits_{\d\downarrow 0}\sup\limits_{n}\mP\left\{2(|a_0|+R)\sup\limits_{s\leq t\leq s+\d}\left|\int_{s}^{t}g_n(u,X^{(n)}_u)\dif W(u)\right|\geq \t^2/3\right\}=0.
$$
By the same deduction to that for $\(\int_s^tg_n(u,X^{(n)}_u)\dif W(u)\)_{s\leq t\leq s+\d}$, one can obtain that 
$$
\lim\limits_{\d\downarrow 0}\sup\limits_{n}\mP\left\{2\sup\limits_{s\leq t\leq s+\d}\left|\int_{s}^{t}\< X^{(n)}(u),g_n(u,X^{(n)}_u)\dif W(u)\>\right|\geq \t^2/3\right\}=0.
$$
Letting $\d\rightarrow 0$ first and then $R\rightarrow\infty$, we infer that
\ce
\lim\limits_{\d\downarrow 0}\sup\limits_{n}\mP\left\{\sup\limits_{0\leq s\leq t\leq s+\d\leq T}|X^{(n)}(t)-X^{(n)}(s)|\geq\t\right\}=0.
\de
Besides, it holds that
\ce
&&\mP\left\{\sup\limits_{-r_0\leq s\leq t\leq s+\d\leq T}|X^{(n)}(t)-X^{(n)}(s)|\geq\t\right\}\\
&\leq&\mP\left\{\sup\limits_{-r_0\leq s\leq t\leq s+\d\leq 0}|X^{(n)}(t)-X^{(n)}(s)|\geq\t/2\right\}\\
&&+\mP\left\{\sup\limits_{0\leq s\leq t\leq s+\d\leq T}|X^{(n)}(t)-X^{(n)}(s)|\geq\t/2\right\}\\
&=&\mP\left\{\sup\limits_{-r_0\leq s\leq t\leq s+\d\leq 0}|\xi(t)-\xi(s)|\geq\t/2\right\}\\
&&+\mP\left\{\sup\limits_{0\leq s\leq t\leq s+\d\leq T}|X^{(n)}(t)-X^{(n)}(s)|\geq\t/2\right\}.
\de
So, the above deduction and the continuity of $\xi$ in $t$ imply (\ref{pxeuxes}). 

Note that for $t\in[0,T]$
\ce
K^{(n)}(t)=\xi(0)+\int_0^t f_n(s,X^{(n)}_s)\dif s+\int_0^tg_n(s,X^{(n)}_s)\dif W(s)-X^{(n)}(t).
\de
Thus, by (\ref{pxeuxes}) we conclude (\ref{pkeukes}).

Combining the above deduction, we obtain that $\{(X^{(n)}_\cdot, K^{(n)}(\cdot)), n\in\mN\}$ is tight in $C([0,T],\sC)\times C([0,T],\mR^d)$.

{\bf Step 3.} We establish the existence of weak solutions for Eq.(\ref{eq1}).

Set
\ce
&&\hat{\Omega}:=C([0,T],\sC)\times C([0,T],\mR^d),\\
&&\hat{\sF}:=\sB\(C([0,T],\sC)\)\times \sB\(C([0,T],\mR^d)\),\\
&&\hat{\sF}_s:=\sB\(C([0,s],\sC)\)\times \sB\(C([0,s],\mR^d)\), \quad 0\leq s\leq T,\\
&&\hat{\Omega}_0:=\{(x,k)\in\hat{\Omega}: k\in A(x_\cdot(0)), k(0)=0\},
\de
and define the canonical processes $\hat{X}: \hat{\Omega}\rightarrow C([0,T],\sC), \hat{K}: \hat{\Omega}\rightarrow C([0,T],\mR^d)$ by
\ce
\hat{X}(x,k)(\cdot)_\cdot:=x_\cdot, \quad \hat{K}(x,k)(\cdot):=\left\{\begin{array}{l} 
k(\cdot), \quad (x,k)\in \hat{\Omega}_0,\\
0, \quad (x,k)\in \hat{\Omega}\setminus \hat{\Omega}_0.
\end{array}
\right.
\de
By {\bf Step 2.}, there exist a subsequence denoted still by $\{(X^{(n)}_\cdot, K^{(n)}(\cdot))\}$ and a probability measure $\hat{\mP}$ on $\hat{\Omega}$ such that as $n$ tends to $\infty$,
\ce
\mP_n\Longrightarrow \hat{\mP},
\de
where $\mP_n=\mP\circ(X^{(n)}_\cdot, K^{(n)}(\cdot))^{-1}$ and  ``$ \Longrightarrow$" stands for the weak convergence of measures. 

Next, we notice that for any constant $l>0$
\ce
\hat{\mP}(\hat{\Omega}_0)&\geq& \hat{\mP}(\{(x,k)\in\hat{\Omega}_0: |k|_0^T\leq l\})\geq \limsup\limits_{n\rightarrow\infty}\mP_n(\{(x,k)\in\hat{\Omega}_0: |k|_0^T\leq l\})\\
&=&1-\liminf\limits_{n\rightarrow\infty}\mP_n(\{(x,k)\in\hat{\Omega}_0: |k|_0^T>l\})\\
&=&1-\liminf\limits_{n\rightarrow\infty}\mP(\{\omega\in\Omega: (X^{(n)}(\omega)_\cdot,K^{(n)}(\omega)(\cdot))\in\hat{\Omega}_0, |K^{(n)}|_0^T>l\})\\
&\geq&1-\frac{C}{l}.
\de
As $l$ tends to $\infty$, it holds that $\hat{\mP}(\hat{\Omega}_0)=1$. So, let $\cN^{\hat{\mP}}$ be the collection of all $\hat{\mP}$ zero subsets of $\hat{\Omega}$ and $\hat{\sF}^{\hat{\mP}}$ be the $g$-field generated 
by $\hat{\sF}$ and $\cN^{\hat{\mP}}$. Set $\hat{\sF}^{\hat{\mP}}_s:=\cap_{v>s}(\hat{\sF}_v\vee\cN^{\hat{\mP}})$ for $0\leq s\leq T$ and $(\hat{\Omega},\hat{\sF}^{\hat{\mP}},(\hat{\sF}^{\hat{\mP}}_s),\hat{\mP})$ is a complete filtered probability space. 

In the following, since $\mP_n$ is the distribution of $(X^{(n)}_\cdot, K^{(n)}(\cdot))$, we see that
$$
M^{(n)}(t):=\hat X_t(0)-\hat X_0(0)+\hat K(t)-\int_0^t f_n(s, \hat{X}_s) \mathrm{d} s, \quad t \in [0, T]
$$
is a $\mP_n$-martingale with
$$
\langle M_i^{(n)}, M_j^{(n)}\rangle(t)=\sum_{k=1}^m \int_0^t\{(g_n)_{i k}(g_n)_{j k}\}(s, \hat{X}_s) \mathrm{d} s, \quad 1 \leq i, j \leq d.
$$
Since $g_n, f_n$ are bounded uniformly in $n$, and as $n \rightarrow \infty$ converge to $g$ and $f$ locally uniformly, by the weak convergence of $\mP_n$ to $\hat{\mP}$ and letting $n \rightarrow \infty$, we see that
$$
M(t):=\hat X_t(0)-\hat X_0(0)+\hat K(t)-\int_0^t f(s, \hat X_s) \mathrm{d} s, \quad t\in[0, T]
$$
is a $\hat{\mP}$-martingale with
$$
\langle M_i, M_j\rangle(t)=\sum_{k=1}^m \int_0^t\{g_{i k} g_{j k}\}(s, \hat X_s) \mathrm{d} s, \quad 1 \leq i, j \leq d.
$$
By \cite[Theorem II.7.1]{iw}, this implies that
$$
M(t)=\int_0^t g(s, \hat X_s) \mathrm{d} \hat{W}(s), \quad t \in[0, T]
$$
for some $m$-dimensional Brownian motion $\hat{W}$ on the filtered probability space $(\hat{\Omega},\hat{\sF}^{\hat{\mP}},(\hat{\sF}^{\hat{\mP}}_s),\hat{\mP})$. Therefore, Eq.(\ref{eq1}) has a weak solution $(\hat{\Omega},\hat{\sF}^{\hat{\mP}},(\hat{\sF}^{\hat{\mP}}_s),\hat{\mP}; \hat W, (\hat X_\cdot, \hat K(\cdot)))$.

{\bf Step 4.} We prove the pathwise uniqueness for Eq.(\ref{eq1}). 

Assume that $(\hat{\Omega},\hat{\sF}^{\hat{\mP}},(\hat{\sF}^{\hat{\mP}}_s),\hat{\mP}; \hat W, (\hat Y_\cdot, \tilde{\hat{K}}(\cdot)))$ is the other weak solution for Eq.(\ref{eq1}) with the same initial value $\hat Y_0=\hat X_0$, that is,
\ce
\hat Y(t)=\hat X_0(0)-\tilde{\hat{K}}(t)+\int_0^tb(s,\hat Y_s)\dif s+\int_0^tg(s,\hat Y_s)\dif W(s).
\de
By the It\^o formula and $({\bf H}'_2)$, we have that
\ce
&&|\hat X(t)-\hat Y(t)|^2 \\
&=&-2\int_0^t\< \hat X(s)-\hat Y(s),\dif \hat K(s)-\dif \tilde{\hat{K}}(s)\> +2\int_0^t\< \hat X(s)-\hat Y(s),f(s,\hat X_s)-f(s,\hat Y_s)\>\dif s \\
&&+2\int_0^t\< \hat X(s)-\hat Y(s),(g(s,\hat X_s)-g(s,\hat Y_s))\dif \hat W(s)\> +\int_0^t\left\|g(s,\hat X_s)-g(s,\hat Y_s)\right\|^2\dif s\\
&\leq&2L'_2\int_{0}^{t}\kappa(\|\hat X_s-\hat Y_s\|_\infty^2)\dif s+2\int_0^t\< \hat X(s)-\hat Y(s),(g(s,\hat X_s)-g(s,\hat Y_s))\dif \hat W(s)\>.
\de
Apply the Burkholder-Davis-Gundy inequality and the Young inequality, we have
\ce
&&\mE\(\sup_{s\in[0,t]}|\hat X(s)-\hat Y(s)|^2\) \\
&\leq&2L'_2\int_{0}^{t}\mE\kappa(\|\hat X_s-\hat Y_s\|_\infty^2)\dif s+12\mE\left(\int_0^t\left|\hat X(s)-\hat Y(s)\right|^2\left\|g(s,\hat X_s)-g(s,\hat Y_s)\right\|^2\dif s\right)^{1/2} \\
&\leq&2L'_2\int_{0}^{t}\mE\kappa(\|\hat X_s-\hat Y_s\|_\infty^2)\dif s+\frac{1}{2}\mE(\sup_{s\in[0,t]}|\hat X(s)-\hat Y(s)|^2) \\ 
&&+C\int_0^t\mE\left(\left\|g(s,\hat X_s)-g(s,\hat Y_s)\right\|^2\right)\dif s\\
&\leq&(2+C)L'_2\int_{0}^{t}\mE\kappa(\|\hat X_s-\hat Y_s\|_\infty^2)\dif s+\frac{1}{2}\mE\left(\sup_{s\in[0,t]}|\hat X(s)-\hat Y(s)|^2\right).
\de 
Furthermore,
\ce
\mE\(\sup_{s\in[0,t]}|\hat X(s)-\hat Y(s)|^{2}\)\leq 2(2+C)L'_2\int_{0}^{t}\mE\kappa(\|\hat X_s-\hat Y_s\|_\infty^2)\dif s.
\de

Besides, on the one hand, since $\hat X(s)=\hat X_0(s)=\hat Y_0(s)=\hat Y(s)$ for $s\in[-r_0,0]$, it holds that
\ce
\mE\(\sup_{s\in[-r_0,t]}|\hat X(s)-\hat Y(s)|^{2}\)=\mE\(\sup_{s\in[0,t]}|\hat X(s)-\hat Y(s)|^{2}\).
\de
On the other hand, by the concavity of $\kappa(\cdot)$, we know that 
\ce
&&\mE\kappa(\|\hat X_s-\hat Y_s\|^2_\infty)\leq \kappa(\mE\|\hat X_s-\hat Y_s\|^2_\infty)\\
&\leq&\kappa(\mE\sup\limits_{u\in[0,s]}\|\hat X_u-\hat Y_u\|^2_\infty)=\kappa\left(\mE\left(\sup\limits_{u\in[-r_0,s]}|\hat X(u)-\hat Y(u)|^2\right)\right). 
\de

Collecting the above deduction, we infer that
\ce
\mE\(\sup_{s\in[-r_0,t]}|\hat X(s)-\hat Y(s)|^{2}\)\leq 2(2+C)L'_2\int_{0}^{t}\kappa\left(\mE\left(\sup\limits_{u\in[-r_0,s]}|\hat X(u)-\hat Y(u)|^2\right)\right)\dif s.
\de
Then by Lemma 3.6 in \cite{dq1}, we get $\mE\(\sup\limits_{s\in[-r_0,t]}\left|\hat X(s)-\hat Y(s)\right|^2\)=0, \forall t\in[0,T]$. That is, $\hat X(t)=\hat Y(t),\forall t\in[-r_0,T]\ a.s.\ \mathbb{P}$.
Finally, by the same deduction to that for $K(t)=\tilde K(t), t\in[0,T]$ in the proof of Theorem \ref{wellpose1}, one can get $\hat K(t)=\tilde{\hat{K}}(t), t\in[0,T]$ under the boundedness and $({\bf H}'_2)$. The proof is complete. 
\end{proof}

At present, we prove Theorem \ref{wellpose2}.

{\bf Proof of Theorem \ref{wellpose2}.} First of all, we take a Lipschitz continuous function sequence $\{h_n\}$ on $\sC$ and a compact set sequence $\{\mK_n\}$ with $\mK_n\uparrow \mR^d$ such that $h_n|_{C([-r_0,0],\mK_n\cap\overline{\cD(A)})}=1, h_n|_{C([-r_0,0],\mK^c_n\cap\overline{\cD(A)})}=0$ for any $n\in\mN$. Set
\ce
f_n(t,\zeta)=h_n(\zeta)f(t,\zeta), \quad g_n(t,\zeta)=h_n(\zeta)g(t,\zeta), \quad t\in[0,T], \quad \zeta\in\sC,
\de
and $f_n, g_n$ are bounded and satisfy $({\bf H}'_2)$. By Proposition \ref{bsboun}, the following path-dependent multivalued SDE
\ce\left\{\begin{array}{ll}
\dif X^{(n)}(t)\in -A(X^{(n)}(t))\dif t+f_n(t,X^{(n)}_t)\dif t+g_n(t,X^{(n)}_t)\dif W(t),\quad t\in[0,T],\\
X^{(n)}_0=\xi \in \sC,
\end{array}
\right.
\de
has a unique strong solution $(X^{(n)}_{\cdot}, K^{(n)}(\cdot))$. Since $h_n=1$ on $\mK_n\cap\overline{\cD(A)}$, the above equation is the same to Eq.(\ref{eq1}) before the solution leaves $\mK_n\cap\overline{\cD(A)}$ and Eq.(\ref{eq1}) has a unique solution $(X_{\cdot}, K(\cdot))$ for $t\in [0, \tau_{\mK_n})$, where $\tau_{\mK_n}:=\inf\{t\in[0,T]; X(t)\notin\mK_n\cap\overline{\cD(A)}\}$. Letting $n\rightarrow\infty$, we conclude that Eq.(\ref{eq1}) has a unique solution $(X_{\cdot}, K(\cdot))$ for $t\in [0,T]$. The proof is complete.

\section{The well-posedness of path-dependent multivalued McKean-Vlasov SDEs}\label{pddexiuni}

In this section, we continue to study the well-posedness of path-dependent multivalued McKean-Vlasov SDEs under the framework of Section \ref{pdexiuni}.

First of all, consider the following path-dependent multivalued McKean-Vlasov SDE:
\be\left\{\begin{array}{ll}
\dif X(t)\in -A(X(t))\dif t+b(t,X_t,\sL_{X_t})\dif t+\s(t,X_t,\sL_{X_t})\dif W(t),\\
X_0=\xi \in \sC,
\end{array}
\right.
\label{eq0}
\ee
where $b:[0,T]\times\sC\times\sP_{2}^{\sC}\to \mR^d$, $ \sigma:[0,T]\times\sC\times\sP_{2}^{\sC}\to \mR^{d\times m}$ are Borel measurable and $\xi$ is a $\sF_0$-measurable $\sC$-valued random variable. Here $X_{\cdot}$ is the segment functional of $X(\cdot)$ defined as $X_t(\theta)=X(t+\theta),\theta\in[-r_0,0]$ and $\sL_{X_t}$ is the distribution of $X_t$. 

Assume:
\begin{enumerate}[$({\bf H}^{\prime\prime}_1)$] 
\item There exists a constant $L^{\prime\prime}_1>0$ such that for any $t\in[0,T]$ and any $\zeta\in\sC$, $\mu\in\sP_2^{\sC}$
$$
|b(t,\zeta,\mu)|^2+\|\s(t,\zeta,\mu)\|^2\leq L^{\prime\prime}_1(1+\|\zeta\|^2_{\infty}+\mu(\|\cdot\|^2_{\infty})).
$$ 
\item For $\forall t\in[0,T]$, $b(t,\cdot,\cdot)$ is continuous on $\sC\times \sP_{2}^{\sC}$, and there exists a constant $L^{\prime\prime}_2>0$ such that for any $t\in[0,T]$ and any $\zeta, \eta\in\sC$, $\mu,\nu\in\sP_2^{\sC}$
\ce
&&2\<b(t,\zeta,\mu)-b(t,\eta,\nu),\zeta(0)-\eta(0)\>\leq L^{\prime\prime}_2\left(\kappa_1\(\|\zeta-\eta\|^2_\infty\)+\kappa_2\(\mW_2(\mu,\nu)^2\)\right),\\
&&\|\s(t,\zeta,\mu)-\s(t,\eta,\nu)\|^2\leq L^{\prime\prime}_2\left(\kappa_1\(\|\zeta-\eta\|^2_\infty\)+\kappa_2\(\mW_2(\mu,\nu)^2\)\right),
\de
where $\kappa_i(x), i=1,2$ are positive, strictly increasing, continuous concave function and satisfy
\ce
\kappa_i(0)=0, \quad \int_{0+}\frac{1}{\kappa_1(x)+\kappa_2(x)}\dif x=\infty.
\de

\end{enumerate}

Now, we state the following theorem, which is the main result in this section.

\bt\label{wellpose3}
Suppose that $({\bf H}^{\prime\prime}_1)$ and $({\bf H}^{\prime\prime}_2)$ hold. Then for any initial value $X_0=\xi$ with $\mE\|\xi\|^2_{\infty}<\infty$, Eq.(\ref{eq0}) has a unique strong solution $(X_{\cdot},K(\cdot))$ and 
$$
\mE\left(\sup_{-r_0 \leq t \leq T}|X(t)|^2\right)<\infty.
$$
\et

In order to prove the above theorem, we construct a sequence of path-dependent multivalued SDEs and study some properties about their solutions. Set
$$
X_t^{(0)}(\theta)=\xi\left(0\wedge(t+\theta)\right) \text{for}\  \theta\in[-r_0,0],\mu_t^{(0)}=\sL_{X_t^{(0)}},0\leq t\leq T.
$$
For any $n\in\mN^+$, consider the following path-dependent multivalued SDE
\be
\left\{\begin{array}{ll}
\dif X^{(n)}(t)\in -A(X^{(n)}(t))\dif t+b(t,X_t^{(n)},\mu_t^{(n-1)})\dif t+\sigma(t,X_t^{(n)},\mu_t^{(n-1)})\dif W(t),\\
X_0^{(n)}=\xi,
\end{array}
\right.
\label{itera}
\ee
where $\mu_t^{(n-1)}:=\sL_{X_t^{(n-1)}}$ and $X_t^{(n)}(\theta):=X^{(n)}(t+\theta)$ for $\theta\in[-r_0,0]$. By Theorem \ref{wellpose2}, we know that Eq.(\ref{itera}) has a unique strong solution denoted as $(X^{(n)}_{\cdot},K^{(n)}(\cdot))$. About $X^{(n)}_{\cdot}$, we have the following result.	

\bl\label{disintpro}
Under $({\bf H}^{\prime\prime}_1)$ and $({\bf H}^{\prime\prime}_2)$, it holds that

$(i)$ for any $n\in\mN^+$
\be
\mE\sup\limits_{s\in[-r_0,t]}|X^{(n)}(s)|^2+\mE|K^{(n)}|_0^t<\infty,\ t\in[0,T];
\label{itera_sec}	
\ee

$(ii)$ $\{X^{(n)}_{\cdot}\}$ is a Cauchy sequence in $L^2(\Omega, \sF, \mP; C([0,T],\sC))$.
\el
\begin{proof}
First of all, we prove $(i)$. For $n=1$ and any $a\in\mathrm{Int}(\cD(A))$, by the It\^o formula, Lemma \ref{L2} and$({\bf H}^{\prime\prime}_1)$, it holds that
\ce
&&\left|X^{(1)}(t)-a\right|^2 \\
&=&\left|\xi(0)-a\right|^{2}-2\int_{0}^{t}\< X^{(1)}(s)-a,\dif K^{(1)}(s)\>+2\int_{0}^{t}\< X^{(1)}(s)-a,b(s,X_{s}^{(1)},\mu_{s}^{(0)})\>\dif s \\
&&+2\int_{0}^{t}\< X^{(1)}(s)-a,\sigma(s,X_{s}^{(1)},\mu_{s}^{(0)})\dif W(s)\>+\int_{0}^{t}\left\|\sigma(s,X_{s}^{(1)},\mu_{s}^{(0)})\right\|^{2}\dif s \\
&\leq&\left|\xi(0)-a\right|^{2}-2\gamma_{1}\left|K^{(1)}\right|_{0}^{t}+2\gamma_{2}\int_{0}^{t}\left|X^{(1)}(s)-a\right|\dif s+2\gamma_{3}t+\int_{0}^{t}\left|X^{(1)}(s)-a\right|^{2}\dif s \\
&&+\int_{0}^{t}\left|b(s,X_{s}^{(1)},\mu_{s}^{(0)})\right|^{2}ds+\int_{0}^{t}\|\sigma(s,X_{s}^{(1)},\mu_{s}^{(0)})\|^{2}ds \\
&&+2\int_{0}^{t}\< X^{(1)}(s)-a,\sigma(s,X_{s}^{(1)},\mu_{s}^{(0)})\dif W(s)\> \\
&\leq&\left|\xi(0)-a\right|^{2}-2\gamma_{1}\left|K^{(1)}\right|_{0}^{t}+2\gamma_{2}\int_{0}^{t}(1+\left|X^{(1)}(s)-a\right|^{2})\dif s+2\gamma_{3}t \\
&&+\int_{0}^{t}\left|X^{(1)}(s)-a\right|^{2}\dif s+L^{\prime\prime}_{1}\int_{0}^{t}\left(1+\left\|X^{(1)}_{s}\right\|_{\infty}^{2}+\mu_{s}^{(0)}(\left\|\cdot\right\|_{\infty}^{2})\right)\dif s \\
&&+2\int_{0}^{t}\< X^{(1)}(s)-a,\sigma(s,X_{s}^{(1)},\mu_{s}^{(0)})\dif W(s)\>.
\de
The Burkholder-Davis-Gundy inequality and the Young inequality imply that
\ce
&&\mE\sup_{s\in[0,t]}|X^{(1)}(s)-a|^{2}+2\gamma_{1}\mE| K^{(1)}|_{0}^{t}\\
&\leq&\mE|\xi(0)-a|^{2}+\(2\gamma_{2}+2\gamma_{3}\)t+(2\g_2+1)\int_{0}^{t}\mE\left|X^{(1)}(s)-a\right|^{2}\dif s\\
&&+L^{\prime\prime}_{1}\int_{0}^{t}\left(1+\mE\left\|X^{(1)}_{s}\right\|_{\infty}^{2}+\mu_{s}^{(0)}(\left\|\cdot\right\|_{\infty}^{2})\right)\dif s\\ 
&&+12\mE\(\int_{0}^{t}|X^{(1)}(s)-a|^{2}\|\sigma(s,X_{s}^{(1)},\mu_{s}^{(0)})\|^{2}ds\)^\frac{1}{2} \\
&\leq&\mE\left|\xi(0)-a\right|^{2}+\(2\gamma_{2}+2\gamma_{3}\)t+(2\g_2+1)\int_{0}^{t}\mE\left|X^{(1)}(s)-a\right|^{2}\dif s\\
&&+L^{\prime\prime}_{1}\int_{0}^{t}\left(1+\mE\left\|X^{(1)}_{s}\right\|_{\infty}^{2}+\mu_{s}^{(0)}(\left\|\cdot\right\|_{\infty}^{2})\right)\dif s\\ 
&&+\frac{1}{2}\mE\(\sup_{s\in[0,t]}|X^{(1)}(s)-a|^{2}\)+C\mE\int_{0}^{t}\|\sigma(s,X_{s}^{(1)},\mu_{s}^{(0)})\|^{2}\dif s.
\de
Moreover, $({\bf H}^{\prime\prime}_1)$ and simple calculation yield that
\be
&&\mE\sup_{s\in[0,t]}|X^{(1)}(s)-a|^{2}+4\gamma_{1}\mE|K^{(1)}|_{0}^{t} \no\\
&\leq& 2\mE\left|\xi(0)-a\right|^2+2(2\gamma_2+2\gamma_3)T+2(2\g_2+1)\int_{0}^{t}\mE\left|X^{(1)}(s)-a\right|^{2}\dif s\no\\
&&+2L^{\prime\prime}_{1}(1+C)\int_{0}^{t}\left(1+\mE\left\|X^{(1)}_{s}\right\|_{\infty}^{2}+\mu_{s}^{(0)}(\left\|\cdot\right\|_{\infty}^{2})\right)\dif s\no\\
&\leq& 2\mE\left|\xi(0)-a\right|^2+2\[2\gamma_2+2\gamma_3+L^{\prime\prime}_1(1+C)\(1+2\left|a\right|^2+\sup_{s\in[0,T]}\mu_s^{(0)}(\|\cdot\|_\infty^2)\)\]T \no\\
&&+2\(2\gamma_2+1+2L^{\prime\prime}_1(1+C)\)\int_0^{t}\mE\sup_{u\in[-r_0,s]}|X^{(1)}(u)-a|^2\dif s.
\label{k=1}
\ee
Note that
\ce
\sup_{s\in[-r_0,t]}|X^{(1)}(s)-a|^{2}&\leq& \sup_{s\in[-r_0,0]}|X^{(1)}(s)-a|^{2}+\sup_{s\in[0,t]}|X^{(1)}(s)-a|^{2}\\
&=& \sup_{s\in[-r_0,0]}|\xi(s)-a|^{2}+\sup_{s\in[0,t]}|X^{(1)}(s)-a|^{2},
\de
and
\ce
\sup_{s\in[0,T]}\mu_s^{(0)}(\|\cdot\|_\infty^2)=\sup_{s\in[0,T]}\mE\|X_s^{(0)}\|_\infty^2=\sup_{s\in[0,T]}\mE\sup_{\t\in[-r_0,0]}|X_s^{(0)}(\t)|^2\leq\mE\|\xi\|_\infty^2.
\de
Thus, it holds that
\be
&&\mE\sup_{s\in[-r_0,t]}|X^{(1)}(s)-a|^{2}+4\gamma_{1}\mE|K^{(1)}|_{0}^{t} \no\\
&\leq&3\mE\sup_{s\in[-r_0,0]}\left|\xi(s)-a\right|^2+2\[2\gamma_2+2\gamma_3+L^{\prime\prime}_1(1+C)\(1+2\left|a\right|^2+\mE\|\xi\|_\infty^2\)\]T \no\\
&&+2\(2\gamma_2+1+2L^{\prime\prime}_1(1+C)\)\int_0^{t}\mE\sup_{u\in[-r_0,s]}|X^{(1)}(u)-a|^2\dif s.
\label{k=1k}
\ee
By the Gronwall inequality, we get that
\be
\mE\sup_{s\in[-r_0,t]}\left|X^{(1)}(s)-a\right|^2\leq C'_1\exp(C'_2t),
\label{x1es}
\ee
where 
\ce
C'_1&=&6\mE\|\xi\|_\infty^2+6|a|^2+2\[2\gamma_2+2\gamma_3+L^{\prime\prime}_1(1+C)\(1+2\left|a\right|^2+\mE\|\xi\|_\infty^2\)\]T,\\
C'_2&=&2\(2\gamma_2+1+2L^{\prime\prime}_1(1+C)\).
\de 
Inserting (\ref{x1es}) into (\ref{k=1k}), we conclude that
\ce
\mE\sup_{s\in[-r_0,t]}|X^{(1)}(s)-a|^{2}+\mE|K^{(1)}|_{0}^{t}<\infty.
\de
	
Now, assume that (\ref{itera_sec}) holds for $n=k$ for some $k\geq 1$. Then we prove that it holds for $n=k+1$. By the similar derivation to that for (\ref{k=1}), for any $a\in\mathrm{Int}(\cD(A))$, we have that
\ce
&&\mE\sup_{s\in[0,t]}|X^{(k+1)}(s)-a|^{2}+4\gamma_{1}\mE|K^{(k+1)}|_{0}^{t} \no\\
&\leq& 2\mE\left|\xi(0)-a\right|^2+2\[2\gamma_2+2\gamma_3+L^{\prime\prime}_1(1+C)\(1+2\left|a\right|^2+\sup_{s\in[0,T]}\mu_s^{(k)}(\|\cdot\|_\infty^2)\)\]T \no\\
&&+2\(2\gamma_2+1+2L^{\prime\prime}_1(1+C)\)\int_0^{t}\mE\sup_{u\in[-r_0,s]}|X^{(k+1)}(u)-a|^2\dif s.
\de
Besides, it holds that 
\ce
\sup_{s\in[-r_0,t]}|X^{(k+1)}(s)-a|^{2}&\leq& \sup_{s\in[-r_0,0]}|X^{(k+1)}(s)-a|^{2}+\sup_{s\in[0,t]}|X^{(k+1)}(s)-a|^{2}\\
&=& \sup_{s\in[-r_0,0]}|\xi(s)-a|^{2}+\sup_{s\in[0,t]}|X^{(k+1)}(s)-a|^{2},
\de
and
\ce
\sup_{s\in[0,T]}\mu_s^{(k)}(\|\cdot\|_\infty^2)=\sup_{s\in[0,T]}\mE\|X_s^{(k)}\|_\infty^2=\sup_{s\in[0,T]}\mE\sup_{\t\in[-r_0,0]}|X_s^{(k)}(\t)|^2.
\de
So, the Gronwall inequality yields that
$$
\mE\sup_{s\in[-r_0,t]}\left|X^{(k+1)}(s)\right|^2+\mE|K^{(k+1)}|_{0}^{t}<\infty.
$$ 
So, we have proved (\ref{itera_sec}) holds for any $n\in\mN^+$.
	
Next, we deal with $(ii)$. By the It\^o formula, Lemma \ref{L1} and $({\bf H}^{\prime\prime}_2)$, it holds that for any $n, m\in\mN^+$, 
\ce
&&|X^{(n+1)}(t)-X^{(m+1)}(t)|^2 \\
&=&-2\int_0^t\< X^{(n+1)}(s)-X^{(m+1)}(s),\dif K^{(n+1)}(s)-\dif K^{(m+1)}(s)\> \\
&&+2\int_0^t\< X^{(n+1)}(s)-X^{(m+1)}(s),b(s,X_s^{(n+1)},\mu_s^{(n)})-b(s,X_s^{(m+1)},\mu_s^{(m)})\>\dif s \\
&&+2\int_0^t\< X^{(n+1)}(s)-X^{(m+1)}(s),\left(\sigma(s,X_s^{(n+1)},\mu_s^{(n)})-\sigma(s,X_s^{(m+1)},\mu_s^{(m)})\right)\dif W(s)\> \\
&&+\int_0^t\|\sigma(s,X_s^{(n+1)},\mu_s^{(n)})-\sigma(s,X_s^{(m+1)},\mu_s^{(m)})\|^2\dif s \\
&\leq&2L^{\prime\prime}_{2}\int_{0}^{t}\left(\kappa_1\left(\left\|X_{s}^{(n+1)}-X_{s}^{(m+1)}\right\|_{\infty}^{2}\right)+\kappa_2\left(\mW_{2}(\mu_{s}^{(n)},\mu_{s}^{(m)})^2\right)\right)\dif s \\
&&+2\int_0^t\< X^{(n+1)}(s)-X^{(m+1)}(s),\left(\sigma(s,X_s^{(n+1)},\mu_s^{(n)})-\sigma(s,X_s^{(m+1)},\mu_s^{(m)})\right)\dif W(s)\>.
\de
By the Burkholder-Davis-Gundy inequality and $({\bf H}^{\prime\prime}_2)$, we obtain that
\be
&&\mE\sup_{s\in[0,t]}\left|X^{(n+1)}(s)-X^{(m+1)}(s)\right|^{2} \no\\
&\leq&2L^{\prime\prime}_{2}\int_{0}^{t}\(\kappa_1\left(\mE\|X_{s}^{(n+1)}-X_{s}^{(m+1)}\|_{\infty}^{2}\right)+\kappa_2\left(\mW_{2}(\mu_{s}^{(n)},\mu_{s}^{(m)})^2\right)\)\dif s \no\\
&&+12\mE\(\int_0^t\left|X^{(n+1)}(s)-X^{(m+1)}(s)\right|^2\left\|\sigma(s,X_s^{(n+1)},\mu_s^{(n)})-\sigma(s,X_s^{(m+1)},\mu_s^{(m)})\right\|^2\dif s\)^{\frac{1}{2}} \no\\
&\leq&2L^{\prime\prime}_{2}\int_{0}^{t}\(\kappa_1\left(\mE\|X_{s}^{(n+1)}-X_{s}^{(m+1)}\|_{\infty}^{2}\right)+\kappa_2\left(\mW_{2}(\mu_{s}^{(n)},\mu_{s}^{(m)})^2\right)\)\dif s\no\\
&&+CL^{\prime\prime}_2\int_0^t\(\kappa_1\left(\mE\|X_{s}^{(n+1)}-X_{s}^{(m+1)}\|_{\infty}^{2}\right)+\kappa_2\left(\mW_{2}(\mu_{s}^{(n)},\mu_{s}^{(m)})^2\right)\)\dif s\no\\
&&+\frac{1}{2}\mE\sup_{s\in[0,t]}\left|X^{(n+1)}(s)-X^{(m+1)}(s)\right|^2 \no\\
&\leq&(CL^{\prime\prime}_2+2L^{\prime\prime}_2)\int_0^t\kappa_1\left(\mE\left\|X_s^{(n+1)}-X_s^{(m+1)}\right\|_\infty^2\right)\dif s\no\\
&&+(CL^{\prime\prime}_2+2L^{\prime\prime}_2)\int_0^t\kappa_2\left(\mE\left\|X_s^{(n)}-X_s^{(m)}\right\|_\infty^2\right)\dif s \no\\
&&+\frac{1}{2}\mE\sup_{s\in[0,t]}\left|X^{(n+1)}(s)-X^{(m+1)}(s)\right|^2,
\label{xn1xmues}
\ee
where in the last inequality we use the fact that $\mW_2(\mu_s^{(n)},\mu_s^{(m)})^2\leq \mE\|X_{s}^{(n)}-X_{s}^{(m)}\|_{\infty}^{2}$. Moreover, it holds that
\ce
&&\mE\sup_{s\in[0,t]}\left|X^{(n+1)}(s)-X^{(m+1)}(s)\right|^{2}\\
&\leq& 2(CL^{\prime\prime}_2+2L^{\prime\prime}_2)\int_0^t\kappa_1\left(\mE\sup\limits_{u\in[-r_0,s]}|X^{(n+1)}(u)-X^{(m+1)}(u)|^2\right)\dif s\\
&&+2(CL^{\prime\prime}_2+2L^{\prime\prime}_2)\int_0^t\kappa_2\left(\mE\sup\limits_{u\in[-r_0,s]}|X^{(n)}(u)-X^{(m)}(u)|^2\right)\dif s.
\de
Note that $X^{(n+1)}(s)=X^{(m+1)}(s)=\xi(s)$ for $s\in[-r_0,0]$. Thus
\ce
\sup_{s\in[-r_0,t]}\left|X^{(n+1)}(s)-X^{(m+1)}(s)\right|^{2}=\sup_{s\in[0,t]}\left|X^{(n+1)}(s)-X^{(m+1)}(s)\right|^{2},
\de
and
\ce
&&\mE\sup_{s\in[-r_0,t]}\left|X^{(n+1)}(s)-X^{(m+1)}(s)\right|^{2}\\
&\leq& 2(CL^{\prime\prime}_2+2L^{\prime\prime}_2)\int_0^t\kappa_1\left(\mE\sup\limits_{u\in[-r_0,s]}|X^{(n+1)}(u)-X^{(m+1)}(u)|^2\right)\dif s\\
&&+2(CL^{\prime\prime}_2+2L^{\prime\prime}_2)\int_0^t\kappa_2\left(\mE\sup\limits_{u\in[-r_0,s]}|X^{(n)}(u)-X^{(m)}(u)|^2\right)\dif s.
\de
Set
\ce
\cG(t):=\varlimsup\limits_{n,m\rightarrow\infty}\mE\sup_{s\in[-r_0,t]}\left|X^{(n)}(s)-X^{(m)}(s)\right|^{2},
\de
and we have that
\ce
\cG(T)\leq 2(CL^{\prime\prime}_2+2L^{\prime\prime}_2)\int_0^T\left(\kappa_1(\cG(s))+\kappa_2(\cG(s))\right)\dif s,
\de
which together with \cite[Lemma 2.1]{zhangx1} yields that $\cG(T)=0$. So, we infer that
\ce
\lim\limits_{n,m\rightarrow\infty}\mE\sup_{s\in[-r_0,T]}|X^{(n)}(s)-X^{(m)}(s)|^{2}=0,
\de
which implies that $\{X^{(n)}_{\cdot}\}$ is a Cauchy sequence in $L^2(\Omega, \sF, \mP; C([0,T],\sC))$. The proof is complete.

Now, it is the position to prove Theorem \ref{wellpose3}.

{\bf Proof of Theorem \ref{wellpose3}.} First of all, we prove the existence. By $(ii)$ in Lemma \ref{disintpro}, there exists a process $(X_t)_{t\in[0,T]}\in L^2(\Omega, \sF, \mP; C([0,T],\sC))$ satisfying 
$$
\lim\limits_{n\to\infty}\mE\sup\limits_{t\in[0,T]}\left\|X^{(n)}_t-X_t\right\|^2_\infty=0.
$$
Since
\ce
\mathbb{W}_2(\mu_t^{(n)},\mu_t)^2\leq \mE\left\|X^{(n)}_t-X_t\right\|^2_\infty,
\de
where $\mu_t$ is the distribution of $X_t$, we know that
$$
\lim_{n\to\infty}\sup_{t\in[0,T]}\mathbb{W}_2(\mu_t^{(n)},\mu_t)^2=0.
$$
 Then, $({\bf H}^{\prime\prime}_1)$,$({\bf H}^{\prime\prime}_2)$ and the dominated convergence theorem imply that $\forall t\in[0,T]$
\ce
\int_{0}^{t}b(s,X^{(n)}_s,\mu^{(n-1)}_s)\dif s &\to&\int_{0}^{t}b(s,X_s,\mu_s)\dif s,\\
\int_{0}^{t}\sigma(s,X^{(n)}_s,\mu^{(n-1)}_s)\dif W(s)&\to& \int_{0}^{t}\sigma(s,X_s,\mu_s)\dif W(s).
\de
Set 
$$
K(t):=\xi(0)-X(t)+\int_0^tb(s,X_s,\mu_s)\dif s+\int_0^t\sigma(s,X_s,\mu_s)\dif W(s),
$$
and by Lemma \ref{L3} and (\ref{itera_sec}), we know that $(X(\cdot),K(\cdot))\in\sA$. Therefore, $(X,K)$ solves Eq.(\ref{eq0}) with 
$$
\mE\sup_{s\in[-r_{0},t]}|X^{(n)}(s)|^{2}<\infty,\forall t\in[0,T].
$$
	
Finally, we prove the pathwise uniqueness. Assume that $(\Omega, \sF, \{\sF_t\}_{t\in[0,T]}, \mP; W, (Z_\cdot,\hat{K}(\cdot)))$ is the other weak solution for Eq.(\ref{eq1}) with the same initial value $Z_0=\xi$, that is
\ce
Z(t)=\xi(0)-\hat{K}(t) +\int_{0}^{t}b(s, Z_s, \sL_{Z_s}) \dif s +\int_{0}^{t} \sigma(s, Z_s, \sL_{Z_s}) \dif W(s).
\de
By the similar deduction to that for (\ref{xn1xmues}), it holds that
\ce
\mE\sup_{s\in[-r_{0},t]}\left|X(s)-Z(s)\right|^{2}&\leq& 2(CL^{\prime\prime}_2+2L^{\prime\prime}_2)\int_0^t\kappa_1(\mE\sup_{u\in[-r_{0},s]}\left|X(u)-Z(u)\right|^{2})\dif s\\
&&+2(CL^{\prime\prime}_2+2L^{\prime\prime}_2)\int_0^t\kappa_2(\mE\sup_{u\in[-r_{0},s]}\left|X(u)-Z(u)\right|^{2})\dif s
\de
And by \cite[Lemma 2.1]{zhangx1}, we have that
$$
\mE\sup_{s\in[-r_{0},t]}\left|X(s)-Z(s)\right|^{2}=0,\forall t\in[0,T].
$$
That is, $X(t)=Y(t),\forall t\in[-r_0,T]\ a.s.\ \mathbb{P}$.

For $\forall t\in[0,T]$,
\ce
K(t) &=&\xi(0)-X(t)+\int_{0}^{t}b(s,X_{s},\sL_{X_s})\dif s+\int_{0}^{t}\sigma(s,X_{s},\sL_{X_s})\dif W(s)\no \\
&=&\xi(0)-Y(t)+\int_{0}^{t}b(s,Z_{s},\sL_{Z_s})\dif s+\int_{0}^{t}\sigma(s,Z_{s},\sL_{Z_s})\dif W(s)\no \\
&=&\hat{K}(t).
\de
Moreover, $K(t), \hat{K}(t)$ are continuous in $t$. Thus, $K(t)=\hat{K}(t),\forall t\in[0,T]\  a.s.\mathbb{P}$. The proof is complete.
\end{proof}


\bigskip
\begin{thebibliography}{999}

\bibitem{Cepa95} E. C\'epa: \'Equations diff\'erentielles stochastiques multivoques, {\it S\'eminaire de probabilit\'es de Strasbourg}, 29(1995)86-107.

\bibitem{Cepa98} E. C\'epa: Probl\'eme de Skorohod multivoque, {\it Ann. Probab.}, 26(1998)500-532. 

\bibitem{chi}  H. Chi: Multivalued stochastic Mckean-Vlasov equation, {\it Acta Math. Sci.}, 34(2014)1731-1740.

\bibitem{dq1} X. Ding and H. Qiao: Euler-Maruyama approximations for stochastic McKean-Vlasov equations with non-Lipschitz coefficients, {\it Journal of Theoretical Probability}, 34(2021)1408-1425.

\bibitem{dq2} X. Ding and H. Qiao: Stability for stochastic McKean-Vlasov equations with non-Lipschitz coefficients, {\it SIAM J. Control Optim.}, 59(2021)887-905.

\bibitem{flqz} K. Fang, W. Liu, H. Qiao and F. Zhu: Asymptotic behaviors of small perturbation for multivalued McKean-Vlasov stochastic differential equations, {\it Applied Mathematics and Optimization}, 88(2023)22.

\bibitem{huangx1} X. Huang: Path-distribution dependent SDEs with singular coefficients, {\it Electron. J. Probab.}, 26(2021)1-21.

\bibitem{huangx2} X. Huang: Strong solutions for functional SDEs with singular drift, {\it Stoch. Dyn.}, 18(2018)1850015.

\bibitem{hrenw} X. Huang, P. Ren and F.-Y. Wang: Distribution dependent stochastic differential equations, {\it Frontiers of Mathematics in China}, 16(2021)257-301.

\bibitem{hrockw} X. Huang, M. R\"ockner and F.-Y. Wang: Nonlinear Fokker-Planck equations for probability measures on path space and path-distribution dependent SDEs, {\it Discrete \& Continuous Dynamical Systems-A}, 39(2019)3017-3035. 

\bibitem{iw} N. Ikeda and S. Watanabe: {\it Stochastic Differential Equations and Diffusion Processes,} 2nd ed., North-Holland/Kodanska, Amsterdam/Tokyo, 1989.

\bibitem{lq} M. Liu and H. Qiao: Parameter estimation of path-dependent McKean-Vlasov stochastic differential equations, {\it Acta Mathematica Scientia}, 42(2022)876-886. 

\bibitem{Mao} X. Mao: {\it Stochastic Differential Equations and Applications}, 2nd ed., Woodhead Publishing, Cambridge, 2008.

\bibitem{McKean} H. P. McKean: A class of Markov processes associated with nonlinear parabolic equations, {\it Proceedings of the National Academy of Sciences of the United States of America}, 56(1966)1907-1911. 

\bibitem{Moha} S. E. A. Mohammed: {\it Stochastic Functional Differential Equations}, Pitman, Boston, 1984.

\bibitem{q} H. Qiao: Limit theorems of invariant measures for multivalued McKean-Vlasov stochastic differential equations, {\it Journal of Mathematical Analysis and Applications}, 528(2023)127532.

\bibitem{qg} H. Qiao and J. Gong: Stability for multivalued McKean-Vlasov stochastic differential equations, {\it Front. Math}, 20(2025)90-932.

\bibitem{Mvon} M. von Renesse and M. Scheutzow: Existence and uniqueness of solutions of stochastic functional differential equations, {\it Random Operators and Stochastic Equations}, 18(2010)267-284. 

\bibitem{wangf1} F.-Y. Wang: Distribution dependent SDEs for Landau type equations, {\it Stochastic Processes and their Applications}, 128(2018)595-621.

\bibitem{wangf2} F.-Y. Wang: {\it Harnack Inequalities for Stochastic Partial Differential Equations}, Springer, Berlin, 2013.

\bibitem{wxz} F. Wu, F. Xi and C. Zhu: On a class of McKean-Vlasov stochastic functional differential equations with applications, {\it Journal of Differential Equations}, 371(2023)31-49. 

\bibitem{za} A. Z\u{a}linescu: Weak solutions and optimal control for multivalued stochastic differential equations, {\it Nonlinear Differential Equations and Applications NoDEA}, 15(2008)511-533.

\bibitem{zhangx1} X. Zhang: Euler-Maruyama approximations for SDEs with non-Lipschitz coefficients and applications, {\it J. Math. Anal. Appl.}, 316(2006)447-458.

\bibitem{zhangx2} X. Zhang: Skorohod problem and multivalued stochastic evolution equations in Banach spaces, {\it  Bull. Sci. Math.}, 131(2007)175-217. 

\bibitem{zhaox} X. Zhao: Well-posedness for path-distribution dependent stochastic differential equations with singular drifts, {\it Journal of Theoretical Probability}, 37(2024)3654-3687.

\bibitem{zhengw} W. Zheng: Tightness results for laws of diffusion processes application to stochastic mechanics, {\it Ann. Inst. H. Poincar\'e, Probab. Stat.}, 21(1985)103-124.

\end{thebibliography}
\end{document}